\documentclass[11pt]{article}
\usepackage{amsfonts,amssymb,amsmath}

\usepackage[english]{babel}
\newtheorem{theorem}{Theorem}
\newtheorem{definition}{Definition}

\newtheorem{proof}{Proof}

\newtheorem{proposition}{Proposition}
\newtheorem{lemma}{Lemma}

\newcommand{\beq}{\begin{eqnarray}}
\newcommand{\eeq}{\end{eqnarray}}

\newcommand{\beqt}{\begin{eqnarray*}}
\newcommand{\eeqt}{\end{eqnarray*}}

\newcommand{\be}{\begin{equation}}
\newcommand{\ee}{\end{equation}}

\newcommand{\bl}{\begin{lemma}}
\newcommand{\el}{\end{lemma}}

\newcommand{\bt}{\begin{theorem}}
\newcommand{\et}{\end{theorem}}

\newcommand{\bd}{\begin{definition}}
\newcommand{\ed}{\end{definition}}

\newcommand{\bp}{\begin{proposition}}
\newcommand{\ep}{\end{proposition}}

\newcommand{\bpr}{\begin{proof}}
\newcommand{\epr}{\end{proof}}

\newcommand{\bi}{\begin{itemize}}
\newcommand{\ei}{\end{itemize}}

\newcommand{\ben}{\begin{enumerate}}
\newcommand{\een}{\end{enumerate}}

\newcommand{\Z}{\mathbb Z}

\newcommand{\N}{\mathbb N}

\newcommand{\E}{\mathbb E}

\newcommand{\s}{\ensuremath{\mathcal{S}}}

\newcommand{\om}{\ensuremath{\omega}}
\newcommand{\Om}{\ensuremath{\Omega}}

\newcommand{\la}{\ensuremath{\Lambda}}

\begin{document}
\title{{\bf Parsimonious Description of Generalized Gibbs Measures : Decimation of the 2d-Ising Model}}

\author{Arnaud Le
Ny\footnote{Laboratoire de Math\'ematiques et d'Analyse Appliqu\'ees (LAMA - UMR CNRS 8050), Universit\'e Paris-Est Cr\'eteil (UPEC), 91 avenue du g\'en\'eral de Gaulle, 94010 Cr\'eteil, France. {\em On leave from} Universit\'e de Paris-Sud, Laboratoire de
Math\'ematiques d'Orsay (LMO UMR CNRS 8628). E-mail : arnaud.le-ny@u-pec.fr.}\footnote{Work supported by the ANR project ANR-JCJC-0139-RANDYMECA.}}

\maketitle

\begin{center}
{\bf Abstract}\\
\end{center}
In this paper, we detail and complete the existing characterizations of the decimation of the Ising model on $\Z^2$ in the generalized Gibbs context. We first recall a few features of the Dobrushin program of restoration of Gibbsianness and present the construction of  global specifications consistent with the extremal decimated measures. We use them to consider these renormalized measures as almost Gibbsian measures and to precise its convex set of DLR measures. We also recall the weakly Gibbsian description and complete it using a potential that admits a quenched correlation decay, i.e. a well-defined configuration-dependent length beyond which this potential decays exponentially. We use these results to incorporate these decimated measures in the new framework of parsimonious random fields that has been recently developed to investigate probability aspects related to neurosciences. 

\medskip
\vspace{7cm}



\newpage
\section{Introduction -- Outline of the paper}
{\em Generalized Gibbs measures} have been introduced during the last decades as a partial answer to a program of restoration of Gibbsianness proposed by Dobrushin in 1995 during a talk in Renkum \cite{Dob2}. The motivations of such a program came from the identification by van Enter {\em et al.} \cite{VEFS} of the pathologies of the Renormalization Group (RG)   as the manifestation of possible non-Gibbsianness of the renormalized measure. These RG transformations being mostly regular scaling transformations -- widely used in the theoretical physics of critical phenomena -- and Gibbs measures designed to rigorously describe equilibrium states and phase transitions in mathematical statistical mechanics, this came out as a surprise. More scientifically, once the elements of surprise were dissipated, it was also taken as a hint that 
the Gibbs property, as described in terms of strongly convergent potentials, might be too strong a notion to fully represent equilibrium states according to probabilistic formulations of the laws of thermodynamics. The program sketched by Dobrushin was thus to relax a bit the conditions required to get this Gibbs property, in order to recover families of random fields stable under the scaling transformations of the renormalization group and still regular enough to be rigorously considered as a proper notion of equilibrium states in mathematical statistical mechanics.\\

Among the main goals in equilibrium mathematical statistical mechanics, one aims at deducing macroscopic behaviors from microscopic rules by considering the possible states of the system as probability measures on the space of microscopic configurations, to interpret some of these probability measures as equilibrium states in a probabilistic sense incorporating ideas from thermodynamics, and eventually to recover a proper notion of phase transition thanks to the possibility of describing different macroscopic equilibrium states for the same microscopic (local) behavior. In view of this last goal, one needs thus to get a construction of probability measures on product spaces alternative to the one of Kolmogorov, based on consistent families of marginals that yield uniqueness of the measure when it exists. It also requires to work at infinite volume and led Dobrushin, Lanford and Ruelle \cite{Dob1,LR} to propose a description of  such probability measures using consistent systems of conditional probabilities with respect to the outside of finite sets, fixed under a boundary condition. This description uses a particular type of consistent kernels called {\em specifications} whose vocation is to characterize such  DLR-consistent systems of conditional probabilities, and eventually to try to specify  probability measures in this way. This {\em DLR approach} allows to describe a large class of measures, strictly larger than the class of Gibbs measures for which an extra topological property is needed. This extra property can be either expressed in terms of potentials that should be uniformly absolutely convergent, or in terms of continuity properties of the  prescribed conditional probabilities, in a manner that sees Gibbs measures as a spatial extension of Markov property (called {\em quasilocality}). In \cite{VEFS}, the authors provide a wide catalog of RG-scaling transformations of Gibbs measures for which this continuity property is violated, by exhibiting points at which no version of these conditional probabilities can be made continuous by changing it on zero-measure sets (essential discontinuity). Considering that this failure of quasilocality might be exceptional, the idea around the first part of the Dobrushin program was to relax the definition and to ask converging and/or continuity properties to hold only on full measure sets. In a second step, one should recover proper variational principles in order to consider these new families of measures as equilibrium states.\\

\newpage
 This gave rise to  two main restoration notions called {\em Almost Gibbsianness} and {\em Weak Gibbsianness}, depending (respectively) if one focuses on almost-sure continuity of conditional probabilities or on almost-sure convergence properties of the potential. The latter notion is implied by the former -- also called almost quasilocalilty -- and has been dismissed because of physically irrelevant thermodynamic properties showing that such a weakly Gibbsian property was indeed too weak a notion to represent equilibrium states with full relevance. Nevertheless, its more physical flavor is very handy and an intermediate notion between them has been coined in \cite{VEV} under the terminology {\em intuitively weak Gibbs}. Here, we let this promising notion for further studies and study whether almost Gibbsianness -- shown to be relevant to describe equilibrium states in \cite{FLNR1,KLNR} -- is concerned for the simplest but fundamental transformation for which RG-pathologies have been detected, the decimation of the two-dimensional ferromagnetic Ising model. We describe what is known for the {\em decimated measures} as far as this generalized Gibbs question is concerned, putting together various results on quasilocality \cite{FP}, almost quasilocality \cite{FLNR1} or weak Gibbsianness \cite{MRSM,shl} and trying also to complete the picture in this particular fundamental case.\\

 A few years ago, independently to this Dobrushin program, other families of random fields have been introduced in \cite{CGL,LO} in order to extend Markov properties in another direction by considering spatial counterparts of Variable Length Markov Chains (VLMC). Motivated by the applicative purpose of extending Rissanen's compression algorithm \cite{Ri}, this gave rise to {\em Parsimonious random fields}. These are DLR-measures for which the conditioning depends on a neighborhood that depends itself on the boundary condition. One of the main motivations of this paper is to complete the weakly Gibbsian description of the decimated measures in order to provide a description of this parsimonious type. Using more general methods and some important results of Maes {\em et al.} \cite{MRSM,shl}, we get an a.s. convergent potential with a so-called {\em quenched correlation decay}, i.e. a potential decaying exponentially beyond a length that depends on the boundary condition. Such a length could be used as a range beyond which the spins are weakly dependent, in order to focus first on the conditioning below it  and to use afterwards the exponential decay to control error terms for spins beyond.\\

In Section 2, we describe the DLR approach to model phase transitions on infinite product probability spaces. We define what we consider to be  DLR measures and classify them  through the families of {\em Quasilocal measures, Gibbs measures, Generalized Gibbs measures} and the recently introduced  {\em Parsimonious random fields}. We eventually illustrate these notions by  the Ising and decimated Ising measures on the square lattice $\Z^2$.  In Section 3, we focus on the fundamental concept of {\em Global specifications} and use a general construction of Fern\'andez and Pfister \cite{FP} to build a global specification for the two-dimensional Ising model at any temperature (Theorem \ref{thmglob}). We use it in Section 4 to get a local specification for the decimated Ising model (Theorem \ref{localspe}) and prove that in dimension two the decimated measures -- i.e. the image measures of the extremal phases of the $2d$-Ising model -- are {\em Almost and Weakly Gibbs} (Theorems \ref{AGChoquet} and \ref{thm9}). In Section 5, we detail the weakly Gibbsian result and provide a {\em Parsimonious description} of these renormalized measures, getting a quenched correlation decay for the renormalized potential (Theorem \ref{thm10}). To do so, we introduce for the latter a  quenched (configuration-dependent) length explicitly described using an {\em Amoeba property of the Ising model} coined by Shlosman \cite{shl},  inspired by the work of Maes {\em et al.} \cite{MRSM}.
\section{DLR Description of Probability Measures}

Our configuration space is the usual infinite product probability space of  spin-systems: $$\big(\Omega,\mathcal{F},\rho)=\big(\{-1,+1\}^{\Z^2}, \{-1,+1\}^{\otimes \Z^2},\rho_0^{\otimes \Z^2} \big)$$ where  $\rho_0=\frac{1}{2} \delta_{-1} + \frac{1}{2} \delta_{+1}$ is the {\em a priori} counting measure on $E=\{-1,+1\}$. \\

The sites of $\mathbb{Z}^2$ are generically denoted by $i=(i^{(1)},i^{(2)}), j=(j^{(1)},j^{(2)})$, etc. and ordered lexicographically : $i \leq j$ {\em iff} $i^{(1)} < j^{(1)}$ or $i^{(1)} = j^{(1)}$ and $i^{(2)} < j^{(2)}$. Two sites $i,j \in \Z^2$ are said to be {\em nearest-neighbors} (n.n.)   {\em iff} $||i-j||_1=1$ and a {\em self-avoiding path} -- or simply a {\em path} -- in $\Z^2$ is a finite sequence $\pi=\{i_1,\dots,i_n\}$ of sites such that $i_j$ and $i_k$ are n.n.  {\em iff} $|j-k|=1$.  We call dual of $\Z^2$ the set $\Z^2 + (\frac{1}{2}, \frac{1}{2})$ and define a {\em contour} $\gamma=(r_1,\dots,r_n)$, of length $|\gamma|=n \in \mathbb{N}^*$, to be a sequence of points in the dual such that $(r_j,\dots,r_n,r_1,\dots,r_{j-1})$ is a path for all $j=1,\dots,n$. The lattice sites enclosed by a contour $\gamma$ form its {\em interior} ${\rm Int} \; \gamma \subset \Z^2$ while its complement in $\mathbb{Z}^2$ is its {\em exterior} ${\rm Ext} \; \gamma$. In a family of contours, a particular contour is said to be {\em external} if it is not surrounded by any other contour of the family, while two contours $\gamma$ and $\gamma'$ are said to be {\em mutually external} iff ${\rm Int} \;\gamma \cap {\rm Int} \; \gamma' = \emptyset$.\\


The set of finite, non-empty, subsets of $\Z^2$ is denoted by $\s$ and for any $\Lambda \in \mathcal{S}$,  $(\Omega_\Lambda, \mathcal{F}_\Lambda, \rho_\Lambda)$ is the finite product space $E^\Lambda$ with the usual product structure. The infinite-volume limit $\Lambda \uparrow \s$ is the standard convergence along nets with respect to the set $\s$ directed by inclusion  \cite{Ge,VEFS}. For the sake of simplicity, it will often be restricted to sequences of cubes $(\Lambda_n)_{n \in \N}$, with $\Lambda_n=([-n,n] \times \Z)^2$. The $\sigma$-field $\mathcal{F}_\Lambda$ is also the $\sigma$-algebra generated by the functions $\{\sigma_i, \; i \in \Lambda\}$ and for all $\sigma, \omega \in \Omega$, we denote $\sigma_\Lambda$ and $\omega_\Lambda$ the projections on $\Omega_\Lambda$, and write similarly $\mu_\la$ for the restriction of $\mu \in \mathcal{M}_1^+$ to $(\Omega_\la,\mathcal{F}_\la)$. We also write $\sigma_\Lambda \omega_{\Lambda^c}$ for the configuration which agrees with $\sigma$ in $\Lambda$ and with $\omega$ in $\Lambda^c$. A function $f$ is said to be {\em local} if there exists $\Lambda \in \s$ s.t. $f$ is $\mathcal{F}_\Lambda$-measurable, and we also write it $f \in \mathcal{F}_\Lambda$. We denote by $\mathcal{F}_{\rm loc}$ the set of local functions, whose  uniform closure $\mathcal{F}_{\rm qloc}$ is the set of {\em quasilocal functions} and coincides\footnote{When $\Omega$ is equipped with the product topology, see \cite{Fer,Ge}.}  with the set $C(\Omega)$ of continuous functions on $\Omega$. In our particular spin system set-up, where we denote by $+$ (resp. $-$) the configuration whose value is $+1$ (resp $-1$) at any site, a function $f$ is said to be {\em right-continuous} (resp. {\em left-continuous}) when for each $\omega \in \Omega$, $\lim_{\Lambda \uparrow \mathcal{S}}  f(\omega_\Lambda+_{\Lambda^c})= f(\omega) \; ({\rm resp.} \; \lim_{\Lambda \uparrow \mathcal{S}}  f(\omega_\Lambda-_{\Lambda^c})= f(\omega))$. When $+$ is replaced by a more general configuration $\theta \in \Omega$, one speaks of {\em quasilocality in the direction} $\theta$.\\

We also generically consider infinite subsets $S \subset \Z^2$, for which all the preceding notations defined for finite $\Lambda$ extend naturally ($\Omega_S,\mathcal{F}_S, \rho_S,\sigma_S$, etc.). Important events to be considered are the {\em asymptotic events}, which are the elements of the {\em tail} $\sigma$-{\em algebra} $\mathcal{F}_{\infty}=\cap_{\Lambda \in \mathcal{S}} \mathcal{F}_{\Lambda^c}$. These events typically do not depend on local behaviors and are mostly obtained by some  limiting procedure. We define translations $\{\tau_i: i \in \Z^2\}$ on  configurations $\omega \in \Omega$, defined via $(\tau_i \omega)(j)=\omega_{i+j}$, for all $i,j \in \Z^2$. Its action extends naturally to functions and measures via the actions $\tau_i f (\omega)=f(\tau_i \omega)$ and $\int f d (\tau_i \mu)= \int (\tau_i f) d\mu$. We denote $\mathcal{M}_1^+$ (resp. $\mathcal{M}_{1,\rm{inv}}^+$) the set of (resp. translation-invariant) probability measures on $(\Omega,\mathcal{F})$. We moreover consider a partial order $\leq$ on $\Omega$: $\sigma \leq \omega$ if and only if $\sigma_i \leq \omega_i$ for all $i \in \Z^2$. This order extends to functions: $f$ is called {\em increasing}  when $\sigma \leq \omega$ implies $f(\sigma) \leq f(\omega)$. It induces then stochastic domination on measures and we write $\mu \leq \nu$ if and only if $\mu[f] \leq \nu[f]$ for all $f$  increasing\footnote{We briefly write $\mu[f]$ for the expectation $\E_\mu[f]$ of a measurable function $f$ under a measure $\mu$.}. 

\subsection{DLR Measures, Quasilocality and Gibbs Measures}

In mathematical statistical mechanics, macroscopic states are represented by  measures in  $\mathcal{M}_1^+$. To describe such measures on the infinite product probability space $\Omega$ in a way that would not necessarily lead to a unique probability measure, and then  to mathematically describe phase transitions, one aims  at describing it by the prescription of its conditional probabilities w.r.t.  boundary conditions. In this so-called {\em DLR approach}\footnote{Introduced by Dobrushin \cite{Dob1}, Lanford and Ruelle \cite{LR} in the late 60's, in an alternative way to the Kolmogorov's description in terms of compatible marginals, the latter implying here existence but also uniqueness.}, a candidate to represent such a system of conditional probabilities are families of probability kernels\footnote{Formally introduced by F\"ollmer \cite{Foll0} and Preston \cite{Pr0} in the mid 70's.} called {\em (local) specification}:
\begin{definition}[(Local) Specification]
A {\em specification} $\gamma=\big(\gamma_\Lambda\big)_{\Lambda \in \s}$  on $(\Omega,\mathcal{F})$ is a family of probability kernels  $\gamma_\Lambda : \Omega \times \mathcal{F} \; \longrightarrow \; [0,1];\; (\omega,A) \; \longmapsto \;=\gamma_\Lambda(A \mid \omega)$
s.t. for all $\Lambda \in \mathcal{S}$:
\begin{enumerate}
\item For all $\omega \in \Omega$, $\gamma_\Lambda(\cdot | \omega)$ is a probability measure on $(\Omega,\mathcal{F})$.
\item For all $A \in \mathcal{F}$, $\gamma_\Lambda(A | \cdot)$ is $\mathcal{F}_{\Lambda^c}$-measurable.
\item For all $\omega \in \Omega$, $\gamma_\Lambda(B|\omega)=\mathbf{1}_B(\omega)$ when $B \in \mathcal{F}_{\Lambda^c}$.
\item For all $\Lambda \subset \Lambda' \in \s$, $\gamma_{\Lambda'} \gamma_{\Lambda}=\gamma_{\Lambda'}$ where 
\be \label{DLR0}
\forall A \in \mathcal{F},\; \forall \omega \in \Omega,\;(\gamma_{\Lambda'} \gamma_\Lambda)(A | \omega)=\int_\Omega \gamma_\Lambda(A | \omega') \gamma_{\Lambda'}(d \omega' | \omega).
\ee
\end{enumerate}
\end{definition}
These kernels also acts on functions and on measures: for all $f \in C(\Omega)$ or $\mu \in \mathcal{M}_1^+$,
$$
\gamma_\Lambda f(\omega):=\int_\Omega f(\sigma) \gamma_\Lambda (d \sigma | \omega)=\gamma_\Lambda [f | \omega] \; {\rm and} \; 
\mu \gamma_\Lambda [f] : = \int_\Omega (\gamma_\Lambda f)(\omega) d \mu (\omega)= \int_\Omega \gamma_\Lambda [f | \omega] \mu(d \omega).
$$

A specification is then said to be {\em monotonicity-preserving} or {\em attractive} when, for all $\la \in \s$ and $f$ increasing, $\gamma_\la f$ is an increasing function (of the boundary condition).\\

 Another way\footnote{Taken as definition of consistency by Kozlov \cite{Ko}.}  to express the {\em consistency} condition 4. is the following {\em key-bar displacement property} (\ref{Keybar}) -- proved in e.g. \cite{Fer,ALN2} -- which tells that for non-null specifications the conditioning  can be freely moved as long as one works on ratios of specifications with boundary conditions of numerator and denominator that remain identical. It will be essential to build our almost-surely convergent potentials  in Section 4 and 5. 

\begin{proposition}[DLR Consistency]\label{propkeybar}\cite{Fer,ALN2}
For a given specification $\gamma=(\gamma_\la)_{\la \in \s}$, consider $\Lambda \subset \Delta \in \s$ and configurations $\sigma, \tau, \omega$ such that $\gamma_\Lambda(\sigma_\Lambda | \tau_{\Delta \setminus \Lambda}  \omega_{\Delta^c}) >0$. Then, for any configuration $\tilde{\sigma} \in \Omega$
\be \label{Keybar}
\frac{\gamma_\Delta (\tilde{\sigma}_\Lambda \tau_{\Delta \setminus \Lambda} | \omega_{\Delta^c})}{\gamma_\Delta(\sigma_\Lambda  \tau_{\Delta \setminus \Lambda}  | \omega_{\Delta^c})}= \frac{\gamma_\Lambda(\tilde{\sigma}_\Lambda | \tau_{\Delta \setminus \Lambda} \omega_{\Delta^c})}{\gamma_\Lambda(\sigma_\Lambda  | \tau_{\Delta \setminus \Lambda}   \omega_{\Delta^c})}.
\ee
\end{proposition}

This consistency reminds one stability property by conditionings in a regular systems of conditional probabilities of $\mu \in \mathcal{M}_1^+$, with the important exception that it is here required to hold everywhere and not only almost surely, because the measure is unknown. Different measures can then have their conditional probabilities described by the same specification but on different full measure sets, leaving the door open to a mathematical description of phase transitions, as we shall see below for the ferromagnetic Ising model on the square lattice $\Z^2$.

\begin{definition}[DLR Measures]
A probability measure $\mu$ on $(\Omega,\mathcal{F})$ is said to be consistent with a specification $\gamma$ (or specified by $\gamma$) when for all $A \in \mathcal{F}$ and $\Lambda \in \s$
\be \label{DLR1}
\mu[A|\mathcal{F}_{\Lambda^c}](\omega)=\gamma_\Lambda(A|\omega), \; \mu{\rm -a.e.} \;  \omega \in \Omega.
\ee
Equivalently, $\mu$ is consistent with $\gamma$ if for any $f \in \mathcal{F}_{\rm{loc}}$, for all $\Lambda \in \mathcal{S}$
\be \label{DLR2}
\int (\gamma_\Lambda f) d \mu = \int f d \mu
\ee
or, {\em  in an even shorter form, if and only if $\mu \gamma_\Lambda = \mu, \; \forall \Lambda \in \s$}. We denote by $\mathcal{G}(\gamma)$ the set of measures consistent with $\gamma$. For a translation-invariant specification, $\mathcal{G}_{\rm{inv}}(\gamma)$ is the set of translation-invariant elements of $\mathcal{G}(\gamma)$.
\end{definition}

In our nice\footnote{For {\em standard Borel spaces}, see \cite{Ge} or the constructions of Goldstein \cite{Gold}, Preston \cite{Pr}  and Sokal \cite{Sok}.} framework, for a given probability measure $\mu$, it is always possible to construct a specification $\gamma$ such that $\mu \in \mathcal{G}(\gamma)$, while on the contrary to Kolmogorov's construction based on marginals, existence of a measure for a given specification is not guaranteed\footnote{There exist indeed specifications $\gamma$ for which $\mathcal{G}(\gamma)=\emptyset$, see \cite{Ge,ALN2}.}, nor uniqueness: we shall see that one can also get more than one element, and in such a case we officially say that there is {\em phase transition}. Central in statistical mechanics, this notion is also  determinant to get non-Gibbsianness. One mostly reserves the terminology "phase" to the t.i. extremal elements of $\mathcal{G}(\gamma)$, and "phase transition" indeed still means that there is more than one phase for the given specification.

The following extension of quasilocality to specifications appears to be a good framework to insure the existence of a consistent probability measure\footnote{More precisely, the set $\mathcal{G}(\gamma)$ is even a non-empty compact convex subset of $\mathcal{M}_1^+$ \cite{Ge}.}. {\bf For a given specification $\gamma$, we denote by $\Omega_\gamma$ the sets of {\em points of continuity of $\gamma$} } i.e. the set of {\em good configurations} $\omega$'s such that $\gamma_\Lambda f$ is continuous at $\omega$ for each $\Lambda \in \mathcal{S}$ and any local function $f$. A {\bf quasilocal specification} corresponds to the case where $\Omega_\gamma=\Omega$ :
\begin{definition}[Quasilocal Specifications -- Quasilocal Measures]
A {\em specification $\gamma$ is quasi\-local} if for each $\Lambda \in \s$ and each $f$ (quasi)local, $\gamma_\Lambda f \in C(\Omega)$. A {\em measure $\mu \in \mathcal{M}_1^+$ is quasilocal} if it is consistent with some quasilocal specification. In particular, it always admits a continuous version of its conditional probabilities w.r.t. the outside of any finite set, as a function of the boundary condition $\omega \in \Omega$. 
\end{definition}

This quasilocality property can also be characterized by the absence of points of essential discontinuity, that are discontinuity points such that, when they are taken as boundary conditions, it is not possible to modify the conditional expectations on zero-measure sets to make them continuous and get continuous versions, see \cite{VEFS, Fer, ALN2} for details.

Basic examples of quasilocal measures are {\em Gibbs measures} defined via {\em potentials}\footnote{Potentials model the interaction between particles while Hamiltonians models energies of configurations.}.

\begin{definition}[Potentials] A {\em potential} is a family $\Phi=\big(\Phi_A\big)_{A \in \s}$ s.t. $\forall A \in \s$, $\Phi_A$ is $\mathcal{F}_A$-measurable. It is:
{\em convergent} at  $\omega \in \Omega$ if, for all $\Lambda \in \s$, the following sum is convergent:
$$
\sum_{A \cap \Lambda \neq \emptyset} \Phi_A(\omega) \; < \; + \infty.
$$
It is {\em absolutely convergent} at $\omega$ when $|\Phi|$ is convergent, and {\em uniformly convergent} if this absolute convergence  is uniform in $\omega \in \Omega$. The (tail-measurable) set of points of convergence of a  potential is denoted by $\Omega_\Phi$. A potential is
{\em Uniformly Absolutely convergent} (UAC) iff
\be \label{UAC}
\forall i \in \Z^2,\; \sum_{A \ni i} \sup_\omega |\Phi_A(\omega)| < \infty.
\ee
\end{definition}

The {\em Hamiltonian with free boundary condition} at finite volume $\Lambda \in \s$ is, for all $\sigma \in \Omega$,
$$
H_\Lambda^{\Phi,f}(\sigma) :=\sum_{A \subset \Lambda} \Phi_A(\sigma) 
$$
while, when it is convergent, one defines the {\em Hamiltonian with boundary condition} $\omega \in \Omega$ by
$$
H_\Lambda^\Phi(\sigma | \omega) := \sum_{A \cap \Lambda \neq \emptyset} \Phi_A(\sigma_\Lambda \omega_{\Lambda^c}).
$$
\begin{definition}[Gibbs Specification -- Gibbs Measures]
Let $\Phi$ be a UAC potential. The {\em Gibbs specification} $\gamma^{\beta \Phi}=\big(\gamma_\Lambda^{\beta \Phi}\big)_{\Lambda \in \s}$ with potential $\Phi$, at temperature $\beta^{-1}>0$,  is defined  by
\be \label{BG}
\forall \lambda \in \s,\; \forall \omega \in \Omega, \; \gamma_\Lambda^{\beta \Phi}(d\sigma | \omega)=\frac{1}{Z^{\beta \Phi}_\Lambda(\omega)} \; e^{-\beta H_\Lambda^\Phi(\sigma | \omega)} \rho_\Lambda\otimes \delta_{\omega_{\Lambda^c}} (d \sigma)
\ee
where $Z_\Lambda^{\beta \Phi}(\omega)$ is a normalizing constant called partition function. A measure $\mu$ is a {\em Gibbs measure} if there exists a UAC potential $\Phi$ and $\beta >0$ s.t. $\mu \in \mathcal{G}(\gamma^{\beta \Phi})$.
\end{definition}

 {\bf Gibbs measures are quasilocal and conversely any non-null quasilocal measure can be written in a Gibbsian way}, with a UAC {\em telescoping} potential built\footnote{See also \cite{Fer,Ge,ALN2}.} by Kozlov \cite{Ko}. We shall use this telescoping procedure to describe under a {\em weakly  Gibbsian form} the decimation of the $+$-phase of the Ising model in dimension two in next section. 

\subsection{Generalized Gibbs Measures}

Although the above Gibbs measures  were commonly accepted as the relevant notion to describe equilibrium states  in mathematical statistical mechanics, it is by now clear that the pathologies of the renormalization group detected in the late seventies in theoretical physics  have to be identified as the manifestation of non-Gibbsianness, while they consist in  reasonable and mostly local renormalization transformations of Gibbs measures. It appeared then that this Gibbs property were probably to strong to be a proper characterization of  equilibrium state and  Dobrushin  initiated in 1995 in  Renkum  a {\em Program of restoration of the Gibbs property}, consisting in two parts: Firstly, to enlarge the family of Gibbs measures -- in order to at least include renormalized Gibbs measures -- and secondly to rigorously  interpret them as equilibrium states by a mathematical formulation of the second law of thermodynamics.
The first part of the Dobrushin program  mainly yields two 
restoration notions, depending if one focuses on a relaxation of the
convergence properties of the potential ({\em Weak
Gibbsianness}), or
 on a relaxation of  topological properties ({\em Almost quasilocality} or {\em Almost Gibbsianness}). Recall that the set $\Omega_\gamma$ of {\em good configurations} of
a specification $\gamma$ gathers its points of continuity.

\begin{definition}[Almost Quasilocal (or Almost Gibbs)]
A probability measure $\mu \in \mathcal{M}_1^+$ is said to be {\em almost Gibbs} if there is a specification $\gamma$ s.t. $\mu \in \mathcal{G}(\gamma)$ and $\mu(\Omega_\gamma)=1$.
\end{definition}

The main other, more physical, restoration notion -- weak Gibbsianness -- appeared rapidly to be indeed weaker and eventually not relevant enough to fully describe a satisfactory notion of equilibrium states \cite{KLNR} without any extra topological requirements. As already noticed initially \cite{MRVM} {\bf almost Gibbsianness implies weak Gibbsianness} and a general construction follows the lines of the  Kozlov's potential built from a vacuum potential whose point-wise convergence is in fact insured by continuity in the direction of the vacuum state.

\begin{definition}[Weak Gibbsianness] A probability measure $\nu \in \mathcal{M}_1^+$ is  {\em Weakly Gibbs}
if there exists a potential $\Psi$ and a tail-measurable\footnote{Tail-measurability is required to insure that the partition
function is well-defined.} set 
$\Omega_\Psi$ on which  $\Psi$ is absolutely convergent, of full measure
 $\nu(\Omega_\Psi)=1$, s.t. $\Psi$ is consistent with $\nu$, i.e. $\mu \in \mathcal{G}(\gamma^{\beta \Psi})$ for some $\beta >0$.
\end{definition}

In case of non-quasilocality, essential discontinuities prevent to define the specification from a UAC potential, needed to have a full Gibbsian formalism. The idea of Dobrushin \cite{Dob2}, developed by Maes {\em et al.} \cite{LMV, MV},  is that one may relax the assumptions and first seeks for an effective interaction  well defined on a set of full measure,  to get the validity of the DLR equations (\ref{DLR1})  for the corresponding "Gibbs-like" specification. One way to get it is to consider first a potential with {\em vacuum state}, taken here to be the $+$ configuration: $\Phi$ is said to be a {\em vacuum potential with vacuum state $+ \in \Omega$} when $\Phi_A(\omega)=0$ whenever $\exists i \in  A$ with $\omega_i=+1$. This potential is the basis of weakly Gibbsian restorations -- if one does not insist on absolute convergence -- whereas a point-wise version of the  telescoping procedure of Kozlov is adapted to get absoluteness. Weak Gibbsianness has been proved  for
most of the renormalized measures 
\cite{BKL,LMV,MRSM,MRVM,MV}. Concerning our 2d-decimated measures, Bricmont {\em et al.} used coarse-graining and contours techniques to build such a weakly Gibbs representation\footnote{They prove that there exists disjoint  sets $\Omega_+, \Omega_- \subset \Omega$ such that 
$\nu^+(\Omega_+)=\nu^-(\Omega_-)=1$
and a translation invariant interaction $\Phi$  with $\Omega_{\Phi}=\Omega_+ \cup \Omega_-$ such that $\nu^+$ and $\nu^-$ are weakly Gibbs for the pair $(\Phi,\Omega_\Phi)$. In their framework $\nu^+$ and $\nu^-$ share the same interaction but concentrate on different configurations.}.\\

We follow in Section 4 an approach of Maes {\em et al.} \cite{MRSM}, inspired by Kozlov-Sullivan's techniques \cite{Ko,Su}. We shall briefly outline how generally one gets such potentials, and describe afterwards the telescoping potential for which one even gets a {\em configuration-dependent} decay of correlations \cite{MRSM,MRVM,MV}. This  configuration-dependence reminds us of the following  {\em Parsimonious Random Fields}, recently introduced by L\"ocherbach {\em et al.} \cite{LO,CGL}.

\subsection{Parsimonious Random Fields}

Inspired by Rissanen's {\em Variable Length Markov Chains} (VLMC, \cite{GL,Ri}), L\"ocherbach {\em et al.} consider a new type of random fields where "to predict the symbols within any finite regions, it is necessary to inspect a random number of neighborhood symbols which change with the boundary condition" \cite{LO}. They first call them {\em Variable-neighborhood random fields} or {\em Parsimonious Markov random fields}, providing natural examples falling into the following Markovian-- hence Gibbsian -- framework (Definition 2.4 in \cite{LO}):
\begin{definition}[Parsimonious {\em Markov} Fields \cite{LO}]\label{DefParsMark}
Let $\gamma$ be a specification. One says that  $\mu \in \mathcal{G}(\gamma)$ is a {\em Parsimonious Markov Field} if for any $\Lambda\in \mathcal{S}$ and for $\mu$-a.e. boundary condition $\omega$, there exists a context $C=C(\omega)=C_\Lambda(\omega)$ {\em finite} $\subset \mathbb{Z}^d$ such that
\be \label{Parsimonious}
\gamma_\Lambda(\cdot \mid \omega_{\Lambda^c})=\gamma_\Lambda(\cdot \mid \omega_C)
\ee
and for all $\tilde{C} \subset \mathbb{Z}^d$, if $\gamma_\Lambda(\cdot \mid \omega_{\Lambda^c})=\gamma_\Lambda(\cdot \mid \omega_{\tilde{C}}) $ then $C \subset \tilde{C}$.
\end{definition}

According to this definition, there might be a set of realizations of $\mu$-measure zero so that $|C_\Lambda(\omega)|=\infty$. This is e.g. the case of another example from \cite{CGL}, the so-called {\em Incompletely observed Markov random fields}. As in the latter, we extend Definition \ref{DefParsMark} allowing the context to be infinite, providing then a new\footnote{While Definition \ref{DefParsMark} is more a new description of particular Gibbs random fields.} family of  possibly non-Gibbsian random fields.

\begin{definition}[Parsimonious {\em Random} Fields \cite{CGL}]\label{DefParsRand}
Let $\gamma$ be a specification. One says that  $\mu \in \mathcal{G}(\gamma)$ is a {\em Parsimonious Random Field} if for any $\Lambda\in \mathcal{S}$ and for $\mu$-a.e. boundary condition $\omega$, there exists a context $C=C(\omega)=C_\Lambda(\omega)$  $\subset \mathbb{Z}^d$ such that
\be \label{Parsimonious2}
\gamma_\Lambda(\cdot \mid \omega_{\Lambda^c})=\gamma_\Lambda(\cdot \mid \omega_C)
\ee
and for all $\tilde{C} \subset \mathbb{Z}^d$, if $\gamma_\Lambda(\cdot \mid \omega_{\Lambda^c})=\gamma_\Lambda(\cdot \mid \omega_{\tilde{C}}) $ then $C \subset \tilde{C}$.
\end{definition}

In \cite{CGL},  the question of finiteness of the interaction neighborhood is related to the absence/presence of phase transition in the underlying Markov field, but in a more intricate way than one could suspect at a first sight. The question of  Gibbsianness or non-Gibbsianness of this model has not yet been  investigated, although we believe that it might be non-Gibbs in the phase transition region by a mechanism similar to the one described in \cite{VEFS} for stochastic RG transformations. In this paper, we investigate whether the decimation of the Ising model could be incorporated in this framework. We shall see in Section 5 that a parsimonious description is indeed possible,  in an  larger class of parsimonious measures.
\subsection{Fundamental Example: Decimation of the $2d$-Ising Model}
\subsubsection{\bf Two-dimensional (ferromagnetic) Ising Model}

We consider thus in this paper the original Ising model, introduced in Ising's thesis directed by Lenz to study ferromagnetism \cite{Ising}. Ising came to the  conclusion of absence of ferromagnetism ordering (i.e. phase transition for us) in dimension 1 but wrongly generalized it to higher dimensions, and the correct result of phase transition in dimension 2 has been settled  a few later after by Peierls. The interaction is {\em nearest-neighbor} (n.n.), with for any $\sigma \in \Omega, \; i,j \in \Z^2$,
\begin{equation}\label{ising}
\Phi_{\{i\}}(\sigma)= - h \; \sigma_i, \; \;  \; \Phi_{\{i,j\}} (\sigma) = - \sigma_i \sigma_j \; {\rm if} \; | i-j | =1
\end{equation}
and $\Phi_A(\sigma)=0$ for other $A \in \mathcal{S}$. The real $h$ is called {\em  magnetic field}.   The existence of a critical temperature  has been qualitatively established by Peierls in 1936 \cite{Peie,Grif} and we state here the results we need through the following theorem  on the structure of the set\footnote{In general, and in particular in our set-up, this set is always a Choquet simplex, i.s. a convex set where each element is determined by a unique convex combination of the extremal elements, see e.g. \cite{Ge, ALN2}. Moreover, these extreme points are the DLR measures that are trivial on the tail $\sigma$-algebra $\mathcal{F}_\infty$.} $\mathcal{G}(\gamma^{\beta \Phi})$ of Gibbs measures for the corresponding   specification $\gamma^{\beta \Phi}$\cite{Aiz,CV,Ge,FP,Hig}.

\begin{theorem}\label{PTIsing} Let $\gamma^{\beta \Phi}$ be the specification (\ref{BG}) with potential (\ref{ising}) at temperature $\beta^{-1}>0$. 
\begin{enumerate}
\item For $h=0$, there exists a critical inverse temperature $0 < \beta_c < + \infty$ such that
\begin{itemize}
\item $\mathcal{G}(\gamma^{\beta \Phi})=\{\mu_\beta \}$ for all $\beta < \beta_c$.
\item $\mathcal{G}(\gamma^{\beta \Phi})=[\mu^-_\beta , \mu^+_\beta]$ for all  $\beta > \beta_c$ where the {\em extremal phases} $\mu^-_\beta \neq \mu^+_\beta$ can be selected via "$-$" or "$+$" boundary conditions: for all $f \in \mathcal{F}_{\rm{qloc}}$,
\be \label{WLimits}
\mu^-_\beta [f]:=\lim_{\Lambda \uparrow \mathcal{S}} \gamma_\Lambda [f \mid -] \; {\rm and} \; \mu^+_\beta [f]:=\lim_{\Lambda \uparrow \mathcal{S}} \gamma_\Lambda [f \mid +].
\ee
\end{itemize}

Moreover, for any $\mu \in \mathcal{G}(\gamma^{\beta \phi})$, for any bounded {\em increasing} $f$, 
\be \label{encadr}
\mu_\beta^-[f] \; \leq \; \mu[f] \; \leq \; \mu_\beta^+[f]
\ee
and the extremal phases have opposite {\em magnetizations} $m^*(\beta):=\mu^+_\beta [\sigma_0]=-\mu^-_\beta [\sigma_0]>0$.
\item For $h>0$ (resp. $h<0$), for any $0<\beta < + \infty$,  $\mathcal{G}(\gamma^{\beta \Phi})=\{\mu_{\beta,h} \}$. The limit $h \to 0$ selects one of the above extreme points: for all $f \in \mathcal{F}_{\rm{qloc}}$, 
$$
\lim_{h \uparrow 0} \mu_{\beta,h}[f ] = \mu^-_\beta [f] \; {\rm and} \; \lim_{h \downarrow 0} \mu_{\beta,h}[f] = \mu^+_\beta [f].
$$
\end{enumerate}

\end{theorem}

Peierls's proof of  {\em ferromagnetic ordering}  at low temperature relies  on a geometric description in terms of {\em contours} that allow to deduce phase transition at low temperature from the behavior at zero temperature. We recall  a few tools that will help to understand the structure of the decimated phase in next section.   The full convex picture has been independently proved by
Aizenmann \cite{Aiz}
 and Higuchi \cite{Hig},  inspired by considerations on
 percolation raised by Russo \cite{Ru}.  An new proof of this  Aizenman-Higuchi theorem, together with a finite-volume version that relies more deeply on the underlying phenomena, has been recently published in \cite{CV}.

Let us first briefly describe the original approach of Peierls who used geometric considerations on contours to get  a temperature-dependent bound on the energy of a configurations that eventually leads to phase transition at low temperature.   Consider $\la \in \s$ to start with the boundary condition $+$  and the probability measure $\gamma_\la^{\beta \Phi}(\cdot| +)$. A {\em contour $\gamma$ of the dual of $\Z^2$ is said to occur} in the configuration $\sigma$, or simply to be a contour {\em of} $\sigma$, if it separates some "$+$" and "$-$" areas of $\omega$, i.e. if $\gamma \subset \big\{ b + (\frac{1}{2},\frac{1}{2}) : b=\{i,j\}, ||i-j||_1=1, \omega_i \neq \omega_j \big\}$. From the Hamiltonian, we see that this event $\gamma$ that a contour occurs requires an energy proportional to its length, so that if $\beta$ is large then long contours will be very improbable w.r.t. the probability $\gamma_\la ^{\beta \Phi}(\cdot | +)$. This yields the following {\em Peierls estimate}:

\be \label{Peierls}
\gamma_\la^{\beta \Phi}(\gamma| +) \leq e^{-2 \beta |\gamma|}.
\ee

From this, it is possible to estimate that the probability that the spin at the origin takes value $-1$, an event which implies the occurrence of contours, by

\be \label{Peierls2}
\gamma_\la^{\beta \Phi}(\sigma_0=-| +) \leq \sum_{l \geq 1} l 3^l e^{-2 \beta l}
\ee
to eventually yield, using (\ref{WLimits}), the weak convergence as $\beta$ goes to infinity of $\mu_\beta^+$  to the Dirac measure $\delta_+$, while the $-$-phase can be similarly proved to converge to the Dirac measure $\delta_-$. One also gets the validity of the Peierls estimate  for either the phases at low enough temperature: the probability that there is a contour which surrounds the origin is bounded by the rhs of (\ref{Peierls2}) while the probability that given contour $\gamma$ surrounds the origin satisfies
\be \label{Peierls1}
\mu_\beta^+(\gamma \ni 0) \leq e^{-c \beta |\gamma|}, \; {\rm for} \; c>0.
\ee
His analysis gave rise to the powerful {\em Pirogov-Sinai theory of phase transitions} for more general models
\cite{PS,Sin},  under a similar {\em Peierls condition}.\\

 In this framework, one gets families of correlated random variables, that are more and more correlated as the temperature decreases and a major challenge is to describe typical properties of configurations that illustrate this picture. The following {\bf Path Large Deviation (PLD) property}, coined by Maes {\em et al.} in \cite{MRSM} is one of these and will be relevant to establish the weak Gibbsian property of the decimated Ising model in the next sections. It precises the convergence of the phases to the above Dirac measures when applied to the magnetization, where for e.g. the $+$-phase, on gets that at low temperature, as $\la \uparrow \s$,
$$
m_\la(\sigma):=\frac{1}{|\Lambda|} \sum_{i \in \la} \sigma_i \; \stackrel{\mu_\beta^+}{\longrightarrow}_\Lambda \; m^*(\beta) >0
$$
so that the event "having the right magnetization" is $\mu_\beta^+$-typical and in particular
\be \label{rightmagn}
\forall \epsilon > 0, \mu_{\la,\beta}^+ \Big[ \big| \frac{1}{|\Lambda|} \sum_{i \in \la} \sigma_i - m^*(\beta) \big| > \epsilon \Big] \longrightarrow_\Lambda \;0.
\ee

Nevertheless, the "$-$"-phase still survives locally -- forming islands of minuses floating inside seas of pluses, see e.g. \cite{DS2} --  and it is no true\footnote{Consider e.g. a path that essentially avoids all contours of the configuration \cite{MRSM}.} that the magnetization remains right along {\em every} path $\pi$. A correct result should involve  long enough paths for typical configurations and Maes {\em et al.} indeed proved the following  {\bf Path Large Deviation property}  (PLD). Denote $\mathcal{A}_{\la}$ the collection of all self-avoiding paths in $\la$, from the origin to the boundary of $\la$. 

\begin{theorem}[PLD]\cite{MRSM, shl}
For $\beta$ large enough and $\epsilon >0$,
\be \label{PLD}
\mu_{\la,\beta}^+ \Big[ \cup_{\pi \in \mathcal{A}_{\la}} \big\{ m^*(\beta) - \frac{1}{|\pi|} \sum_{i \in \pi} \sigma_i > \epsilon \big\} \Big] \; \longrightarrow_\la \; 0.
\ee
\end{theorem}
 The typical configurations where it holds, said to have the {\em  Path Large Deviation property}, will be the source of good configurations for the decimation of the $+$-phase. Satisfying this property means roughly that for {\em almost every} configuration, {\em every} site $i$ and {\em every} finite self-avoiding path $\pi$ starting from $i$, the average magnetization along $\pi$ is not to low. In fact,  \cite{shl} requires it to be valid only for paths whose length $|\pi|$ exceeds a certain threshold value $l_i(\sigma)$ and to express this length, we use in next section the {\em Amoeba property}, weaker than PLD but which also holds typically. It will be interpreted as a length to characterize the quenched correlation decay of the renormalized potential in Section 5.  \\

The contour estimate (\ref{Peierls1}), combined with percolation techniques, allows  to describe more deeply typical properties of the measure, and correlations under it. Using also coupling  techniques, Burton {\em et al.} \cite{BS} proved  that at low enough temperature\footnote{A stronger statement is true at high temperatures and it is conjectured to be true at any non-critical temperature, while it is known that at $\beta_c$ the correlations decay sub-exponentially according to a power law.}, $\mu:=\mu_\beta^+$ is {\bf Quite Weak Bernoulli with Exponential rate} (QWBE): under the total variation norm $|| \cdot ||_{{\rm TV}}$,
\be \label{QWBE}
\forall n  \geq 0, \forall \epsilon > 0, \exists \lambda_\epsilon >0, C_\epsilon > 1, \Big| \Big| \mu_{\la_n^c \cup \la_{[n(1-\epsilon)]}} - \mu_{\la_n^c} \times \mu_{\la_{[n(1-\epsilon)]}} \Big| \Big|_{{\rm TV}} \leq c_\epsilon e^{- \lambda_\epsilon n}.
\ee
We shall adapt their techniques in Section 5 to get a {\bf Quenched Correlation Decay} (QCD) of constrained measures, following Maes {\em et al.} \cite{MRSM,shl}.

\subsubsection{\bf Decimation of the 2d-Ising model}

This basic example already captures the
main non-trivial features of the pathologies of the renormalization group (RG) and can also be seen as the
projection of the 2d-Ising model on the sub-lattice of even sites
\cite{VEFS}, defined as 
\be \label{DefDec}
T \colon (\Omega,\mathcal{F}) \longrightarrow (\Omega',\mathcal{F}')=(\Omega,\mathcal{F}); \; 
\omega \; \;   \longmapsto \omega'=(\omega'_i)_{i \in
\mathbb{Z}^2}, \; {\rm with} \;  \omega'_{i}=\omega_{2i}
\ee

\begin{theorem}\cite{VEFS}\label{Dec} Consider the  Ising model (\ref{ising}) on  $\Z^2$ with zero magnetic field at inverse temperature $\beta  > \tilde{\beta}_c=\frac{1}{2} \cosh^{-1}\big(e^{2 \beta_c}\big)$ 
and denote by $\nu_\beta := T \mu_\beta$ the decimation of any
Gibbs measure $\mu_\beta$ for this
model. Then $\nu_\beta$ is not
quasilocal, hence non-Gibbs.
\end{theorem}

For any $\om \in \Omega$, we denote thus $\om' \in \Om'=\{\pm 1 \}^{\Z^2}$ its image by this decimation transformation, i.e. its restriction on the even lattice rescaled afterwards, while one denotes $\xi=\xi(\om) \in \Omega_2:=\{\pm 1\}^{2 \Z^2}$ its restriction before rescaling.  To prove non-quasilocality of the renormalized measure $\nu_\beta$, one exhibits a {\em bad configuration} where the conditional expectation  of a local function,  w.r.t. the outside of a finite set, is {\em essentially discontinuous}, i.e. discontinuous on a neighborhood\footnote{The former is  automatically of non-zero $\nu^+_\beta$-measure in our topological settings  and thus modifying conditional probabilities on a negligible set cannot make it continuous.}. The role of a bad configuration is played here by the {\em alternating configuration} $\omega'_{\rm{alt}}$ defined for all $i=(i_1,i_2) \in \mathbb{Z}^2$ by $(\omega'_{\rm alt})_i=(-1)^{i_1+i_2}$. Computing the magnetization under the decimated measure, conditioned on $\omega'_{{\rm alt}}$ outside the origin, gives different limits when one approaches it with all $+$ (resp. all $-$) arbitrarily far away, as soon as phase transition is possible in the decorated lattice, the planar lattice got from $\mathbb{Z}^2$ by removing the even sites. The critical temperature $\tilde{\beta}_c$ on this lattice is known explicitly \cite{Sy}, and one recovers thus non-Gibbsianness  for $\beta > \tilde{\beta}_c$.\\

Thus, the very peculiar -- and atypical -- alternating configuration is a point of discontinuity for the image measure because the measure obtained by constraining the $+$-phase to be in this configuration on the even sites exhibits a phase transition: a possible long range order in the "internal" lattice can carry dependencies from arbitrarily far away. The main question is now to find other discontinuity points and to investigate how (a)typical they are. To get non-quasilocality at an image configuration $\omega' \in \Omega'$ seems to be intimately related to the occurrence of phase transition of the measure $\mu^\xi$, constrained to coincide with $\xi=\xi(\omega')$ on the even sites\footnote{One sometimes speaks about a {\em constrained} or {\em hidden} phase transitions.}. The global neutrality of the alternating configuration allowed it, and it seems important that the following\footnote{This set has been originally introduced in \cite{BS} to get  correlation estimates for the $+$-phase.} {\em sparseness set},  composed of the islands of $-1$'s, 
\be \label{spareset}
D(\xi)= \big\{ j \in 2\Z^2: \xi_j=-1 \big\} 
\ee
 is not too sparse. On the other side, by conditioning with the all $+$ configuration on the even sites, one gets a constrained model that does not allow phase transition\footnote{The conditioning acts as an  external field $h_{2i}=+1$ on the even sites, known to lead to uniqueness with the $+$-phase as an equilibrium state.} and  similarly, one can easily see that a configuration with only a finite number of minuses is good.

 A relevant condition has been coined by Shlosman \cite{shl} in terms of  contours. Called {\em Amoeba property}, it has been designed to characterize configurations $\xi$ such that for every long enough path $\pi$, the ratio $|\pi \cap D(\xi)|/|\pi|$ is small and, for every $\la \in \s$ large enough, with $|\la \cap 2\Z^2 |$ of the order of $\la$, the ratio $|\la \cap D(\xi)|/|\la \cap 2\Z^2|$ is also small. First, we introduce an {\em amoeba} consisting in an external contour $\Gamma$, candidate to enclose all the minuses of a configuration, and internal ones that should enclose all the islands of plusses:
 \begin{definition}\label{Amoeba}
 A collection $G$ of compatible contours $\{\Gamma, (\gamma_k)_{k=1,\dots,r}\}$ is called an {\em Amoeba} iff
 \begin{enumerate}
 \item The {\em contour} $\Gamma$ is a unique exterior contour of the collection $G$.
 \item Every  $\gamma_k$ surrounds at least one site belonging to $2 \Z^2$ and the $ \gamma_k$'s  are mutually external.
\end{enumerate}
\end{definition}
Amoebas are used to characterize the spareness of $D(\xi)$. We first define compatible amoebas and to do so denote by $D^+(\xi)$ its complement in $2 \Z^2$.
\begin{definition}
Let $\xi=(\xi_j)_{j \in 2 \Z^2} \in \Omega_2:=\{\pm 1\}^{2 \Z^2}$ with spareness set $D(\xi)$.  An amoeba $G=(\Gamma,(\gamma_k)_{k=1,\dots,r})$ is said to be {\em compatible with the configuration} $\xi$  if 
\begin{enumerate}
\item The  contours of $G$ are contours of a configuration $\sigma \in \Omega$ compatible with $\xi \in \Omega_2$.
\item $D^+(\xi) \cap {\rm Int} \; \Gamma \subset \cup_{k=1}^r {\rm Int} \; \gamma_k$.
\item The interior of every contour $\gamma_k$ does intersect $2 \Z^2$ at sites where $\xi$ takes value $+1$.
\item Denote $D_\Gamma (\xi)$ the connected component of $(D(\xi) \cap {\rm Int} \Gamma)$ containing the inner boundary of $\Gamma$. Then the interior of every contour does not intersect $D_\Gamma(\xi)$.
\end{enumerate}
\end{definition}
Compatibility implies in particular that $ D(\xi) \cap {\rm Int} G = 2 \Z^2  \cap {\rm Int} G$. 
Write ${\rm diam}(G)={\rm diam}(\Gamma)$, ${\rm Int} \; G= {\rm Int} \; \Gamma \cap \big( \cap_{k=1}^r {\rm Ext} \; \gamma_k \big)$  and $|G|=|\Gamma| + \sum_{k=1}^r |\gamma_k|$. Using these definitions, a good way of characterizing spare enough $D(\xi)$ is to ask large enough amoebas compatible with $\xi$ to enclose a much bigger number of sites of $2 \Z^2$ than sites of $D(\xi)$. 
\begin{definition}
An amoeba $G$ compatible with $\xi$ is called {\em benign} iff there exists $\lambda >0$  
\be \label{benign}
| D(\xi) \cap {\rm Int} \; G | \leq \lambda |G|.
\ee
Otherwise it is called {\em malignant}.
\end{definition}
This somehow prevents significant long range orders due to hidden phase transition and allows the definition of a potential when it holds for amoeba large enough. The quenched correlation decay of the transformed potential will be a consequence of the typicality of such behaviors, at least for large enough amoebas, as proved in \cite{shl}.

\begin{definition}\label{AmoebaProp}{\bf (Amoeba Property)}
A configuration $\xi$ satisfies the {\em amoeba property} with the functions $l(\xi)=\big(l_j(\xi)\big)_{j \in \Z^2}$ iff every amoeba $G(\Gamma,\gamma_k)$ compatible with $\xi$, for which $j \in {\rm Int} \;  \Gamma$ and ${\rm diam}(G) \geq l_j(\xi)$, is {\em benign}.
\end{definition}

To get this amoeba property, one proves that "having the right magnetization along every path for $\xi$" implies an amoeba property, and thus that the validity of the PLD (\ref{PLD}) for $\mu_\beta^+$ implies  that amoeba property is $\mu_\beta^+$-typical at low enough temperature. This has been performed by Shlosman in \cite{shl} by extending previous results  got with Maes {\em et al.} \cite{MRSM}.

\begin{theorem}\cite{shl}\label{Amoebathm}
$\forall \lambda >0$, $\exists \beta_0(\lambda) >0$ s.t. for low enough temperature  $\beta \geq \beta_0$, for any $ j \in 2 \Z^2$, any configuration $\xi \in \Omega_{2}$ on the even sublattice of $\Z^2$, there exists configuration-dependent lengths $l_j(\xi) \geq$ and constants $C(\xi) \geq 0$ with
$$l_j(\xi) \leq |j| + C(\xi),\; \; {\rm  such \; that}: $$
\begin{enumerate}
\item The set 
\be \label{tildeOmega}
\tilde{\Omega} = \big\{ \xi \in 2 \Z^2,\; \sup_{j \in 2\Z^2} l_j(\xi) < + \infty \; {\rm  and} \;  C(\xi) < + \infty \big\}
\ee
has a full measure for the decimation of the $+$-phase: $\nu^+_\beta\big(\tilde{\Omega} \big)=1$ for $\beta \geq \beta_0$.
\item Every configuration $\xi \in  \tilde{\Omega}$ satisfies the amoeba property with the function $l_j(\xi)$.
\end{enumerate}
\end{theorem}

We shall describe the lengths $l_j(\xi)$ in the course of the parsimonious description of the decimated measure $\nu^+$ in Section 5. They correspond to the lengths beyond which the amoebas compatible with $\xi$ are benign, and appear indeed to be finite for typical configurations.  This allows afterwards to typically get a Peierls estimate for a constrained measure that yield  a quenched exponential correlation decay by a suitable adaptation of the percolation methods of Burton and Steif \cite{BS}. This nice decay will imply the absolute convergence of the telescoping potential for the typical configurations $\omega$ compatible with the latter typical $\xi's$. We describe it in section 5, before we describe fundamental objects needed to formalize properly our claims and to describe generalized Gibbs measures.\\

 We first need to describe a specification $\gamma^+$ consistent with $\nu^+$.

\section{Global specifications for the 2d-Ising model}

The  proof of non-quasilocality in Theorem \ref{Dec} \cite{VEFS,ALN1} requires  to get conditional probabilities of the original measure w.r.t. to the outside of {\em infinite} sets. To do so, one formally needs {\em global specifications for the original measures}, which extend (local) specifications to infinite sets. Indeed, in the proof sketched above, if one wants to express the conditional probability $\nu[\sigma'_{\bf 0} | \mathcal{F}_{\{{\bf 0}\}^c}]$, where $\bf 0$ is the origin of $\Z^2$, using the image measure $\nu=T \mu$, so that one writes $\nu(A')=\mu(T^{-1}A')$ for any $A' \in \mathcal{F}'$. Thus, one  has to  express $\mu[\sigma_{\bf 0} | \mathcal{F}_{T^{-1}\{{\bf 0}\}}]$ where $T^{-1}\{{\bf 0}\}=2\Z^2 \setminus \{{\bf 0} \}$ is {\em not} the complement of a finite set, but  this of the infinite set $\lambda= (2\Z^2)^c \cup \{{\bf 0}\}$. Thus, we need to express the conditional probability $\mu[\sigma_{\bf 0} | \mathcal{F}_{\lambda^c}]$ with $\lambda$ {\em infinite}, and this is provided by a {\em global} specification with whom $\mu$ is consistent, instead of a local one. 

\begin{definition}[Global specification \cite{FP}]
A {\em Global specification} $\Gamma$ on $\Z^2$ is a family of probability kernels $\Gamma=(\Gamma_S)_{S \subset \Z^2}$ on $(\Omega,\mathcal{F})$ such that for {\em any} $S$ subset of $\Z^2$:
\begin{enumerate}
\item $\Gamma_S(\cdot | \omega)$ is a probability measure on $(\Omega,\mathcal{F})$ for all $\omega \in \Omega$.
\item $\Gamma_S(A | \cdot)$ is $\mathcal{F}_{S^c}$-measurable for all $A \in \mathcal{F}$.
\item{(Properness)}  $\Gamma_S(B|\omega)=\mathbf{1}_B(\omega)$ when $B \in \mathcal{F}_{S^c}$.
\item{(Consistency)}  For all $S_1 \subset S_2 \subset \Z^2$, $\Gamma_{S_2} \Gamma_{S_1}=\Gamma_{S_2}$.
\end{enumerate}
\end{definition}
\begin{definition}
Let $\Gamma$ be a global specification. We write $\mu \in \mathcal{G}(\Gamma)$, or say that  $\mu \in \mathcal{M}_1^+$ is {\em compatible with} $\Gamma$, if for all $A \in \mathcal{F}$ and {\em any} $S \subset \Z^2$,
\be \label{DLR4}
\mu[A|\mathcal{F}_{S^c}](\omega)=\Gamma_S(A|\omega), \; \mu{\rm -a.e.} \;  \omega.
\ee
\end{definition}

The only difference with {\em (local)} specifications is thus that consistency (\ref{DLR4}) has to hold for any subset $S$, including the infinite ones. This implies in particular that there is {\em at most} one  measure $\mu$ compatible with a global specification $\Gamma$, and if there is one, then for all $\omega$ and any $f \in \mathcal{F}_{{\rm loc}}, \; \Gamma_{\Z^2}[f |\omega]=\mu[f | \mathcal{F}_{\emptyset}](\om)=\mu[f]$. On the other hand, starting from a local specification and a consistent measure, Fern\'andez {\em et al.} \cite{FP}  investigate the existence of a global specification such that there is also compatibility at the global level. They provide a positive answer for attractive right-continuous (local) specificatiosn,  by means of a construction inspired from F\"ollmer \cite{Foll} or Goldstein \cite{Gold2}.\\

 We describe it for the $2d$-Ising model, to provide afterwards  local specifications $\gamma^+$ and $\gamma^-$ with whom the decimated measures $\nu^+$ and $\nu^-$ are consistent. We use the construction of \cite{FP} and provide a global specification $\Gamma^+$ (resp $\Gamma^-$) for the original Ising-Gibbs measure\footnote{A similar construction can be done for any $\gamma$  monotonicity-preserving and right- or left-continuous.} $\mu_\beta^+$ (resp. $\mu_\beta^-$) at {\em any} temperature $\beta^{-1} >0$.  We write from now $\gamma^I$ for the local specification $\gamma^{\beta \Phi}$ of the 2d-Ising model, with the n.n. potential $\Phi$ given by (\ref{ising}).

\begin{lemma}\cite{FP}
The Ising specification $\gamma^I$ is right- (resp. left)-continuous and monotonicity-preserving. This allows to construct a global specification $\Gamma^+$ (resp. $\Gamma^-$) such that $\mu_\beta^+ \in \mathcal{G}(\Gamma^+)$ (resp. $\mu_\beta^- \in \mathcal{G}(\Gamma^-)$) . This global specification is monotonicity-preserving and right-continuous (resp. left-continuous), but needs not to be quasilocal in general.
\end{lemma}
When the set $S$ is infinite, one proceeds in two steps, {\bf which order is crucial}: Freeze first the configuration into $\omega$ on $S^c$ and  perform {\em afterwards} the weak limit with $+$- boundary condition {\em in} $S$, to get the constrained measure $\mu^{+,\omega}_S$ on $(\Omega_S,\mathcal{F}_S)$ that defines the kernel corresponding to $S$. This prevents any mandatory quasilocality of the global specification when the corresponding local one is quasilocal, and this indeed sometimes leads to non-Gibbsianness.\\

{\bf Construction of the global specifications} \cite{FP,Foll,Gold2} 

It is a direct consequence of Proposition 3.1 of \cite{FP}, where the full proof is available.  First, one obviously defines $\Gamma^+$ to coincide with the local specification $\gamma^I$ on finite sets, so we put
$$
\forall \Lambda \in \mathcal{S},\; \Gamma^+_\Lambda:=\gamma^I_\Lambda.
$$

To describe its extension on infinite sets, assume first that such a global specification $\Gamma^+$ exists and consider an infinite set $S \subset \Z^2$. By consistency, one has in particular for all $\Lambda$ finite, $ \Gamma_S^+ =  \Gamma_S^+  \Gamma^+_\Lambda =  \Gamma_S^+   \gamma^I_\Lambda  $ so that one should have for any $f \in \mathcal{F}_{\rm{loc}}$ and all $\Lambda \in \mathcal{S}$,
$$
\forall \omega \in \Omega,\; \Gamma_S^+(f |  \omega)=  \int_\Omega \gamma_\Lambda^I(f | \eta) \Gamma_S^+(d \eta | \omega)=\int_\Omega \gamma_\Lambda^I(f | \eta_S\omega_{S^c}) \Gamma_S^+(d \eta | \omega)
$$
where the second equality is got by properness of the kernel $\Gamma_S^+$, which requires that the boundary condition coincides with $\omega$ outside\footnote{One says that the configuration is {\em frozen} in $\omega$ on $S^c$.} $S$. For any $\omega \in \Omega$, $\Gamma^+_S(\cdot \mid \omega)$ has thus to be an infinite volume probability measure that should be consistent with the local specification $\gamma^{S,\omega}$ defined for any configuration $\omega$, for all $\eta$, by
\be \label{frozenspe}
\gamma^{S,\omega}_\Lambda(\cdot \mid \eta):=\gamma_\Lambda^I(\cdot \mid  \eta_S \omega_{S^c}).
\ee

Thus, to get consistency for the global specification, one has to be able to choose for any $\omega$ an element of $\mathcal{G}(\gamma^{S,\omega})$, in a measurable way in order to get a probability kernel, so that in the infinite-volume limit, one recovers consistency with $\mu^+$. In the absence of phase transition, the choice is obvious but we mainly work here in the non-uniqueness region. In this attractive right-continuous case, there is a natural candidate: the weak limit of the kernels restricted, for all $\omega$, at finite volume $\Delta \subset S$, when the outside of $\Delta$ is fixed in the $+$ boundary condition, unless the outside of $S$ that still coincide with $\omega$. There remains  to check consistency, but monotonicity-preservation and right-continuity will do the job. Hence, we select the measure in $\mathcal{G}(\gamma^{S,\omega})$ to be the {\em constrained measure}  on $(\Omega_S,\mathcal{F}_S)$ obtained by weak limit as follows:
\be \label{Constr}
\mu_S^{+,\omega}(\cdot):=\lim_{\Delta \uparrow S} \gamma^I_\Delta (\cdot \mid +_S \omega_{S^c}).
\ee
The weak limit exists by monotonicity, and right-continuity of $\gamma^I$ insures that $\mu_S^{+,\omega} \in \mathcal{G}(\gamma^{S,\omega})$.\\

 Define now for all infinite $S \subset Z^2$ the kernels 
\be \label{GlobalDec}
\Gamma_S^+(d\eta \mid \omega):=\mu_S^{+,\omega} (d \eta_S) \otimes \delta_{\omega_{S^c}}(d \eta_{S^c})
\ee
in such a way that
\be \label{Glob2}
\Gamma_S^+(d\eta | \omega)=\lim_{\Lambda \uparrow S} \gamma_\Lambda^I(d\eta | +_S \omega_{S^c})
\ee
to get a probability measure on $(\Omega,\mathcal{F})$ by construction (Item 1. for the kernel to be a global specification). Measurability w.r.t. $\mathcal{F}_{S^c}$ (Item 2) is insured also by this (monotone) limiting procedure with freezing in $\omega_{S^c}$, which also yields properness (Item 3). It inherits of  monotone preservation from the  original specification and one then proves the right-continuity of this kernel to eventually get the consistency condition for infinite subsets (Item 4: $\mu^+ \in \mathcal{G}(\Gamma^+)$). We stress that the latter requires a careful use of the monotone convergence theorem, see again \cite{FP}. Proceeding similarly starting with the $-$-phase, we get:

\begin{theorem}\cite{FP}\label{thmglob}
Consider the Ising model on $\Z^2$ at inverse temperature $\beta >0$ with specification $\gamma^I$ given by (\ref{BG}) with the n.n. potential (\ref{ising}), and its extremal Gibbs measures $\mu^+$ and $\mu^-$. Define $\Gamma^+=(\Gamma_S^+)_{S \subset \Z^2}$ to be the family of probability kernels on $(\Omega, \mathcal{F})$  as follows:
\begin{itemize}
\item For $S=\Lambda$ finite, for all $\omega \in \Omega$,
$$
\Gamma^+_\Lambda(d \sigma | \omega) := \gamma^I_\Lambda (d \sigma |  \omega)
$$
\item For $S$ infinite, for all $\omega \in \Omega$,
$$
\Gamma^+_S(d\sigma | \omega):=\mu_S^{+,\omega} \otimes \delta_{\omega_{S^c}}(d \sigma)
$$
where the constrained measure $\mu_S^{+,\omega}$  is the weak limit got with freezing in $+_S \omega_{S^c}$ on $\Lambda^c$: $$\mu_S^{+,\omega}(d \sigma_S):=\lim_{\Delta \uparrow S} \gamma^I_\Delta (d \sigma\mid +_S \omega_{S^c}).$$
\end{itemize}
Then $\Gamma^+$ is a global specification such that $\mu^+ \in \mathcal{G}(\Gamma^+)$. It is moreover monotonicity-preserving and right-continuous, but not quasilocal when $\beta > \tilde{\beta}_c >0$. Similarly, one defines a monotonicity-preserving and left-continuous global specification $\Gamma^-$ such that $\mu^- \in \mathcal{G}(\Gamma^-)$.
\end{theorem}

\section{Almost and weak Gibbisanness of the decimated measures}

\subsection{Construction of a local specification consistent with $\nu^+_\beta$ (resp. $\nu^-_\beta$)}

It is based on the global specification $\Gamma^+$ (resp. $\Gamma^-$) using the infinite set $S=2 \mathbb{Z}^2$, or more precisely its complement $(2 \mathbb{Z}^2)^c$ or local modifications of it, as we shall see. Writing DLR equation (\ref{DLR0})  in a short atomic form, one seeks for a (local) specification $\gamma^+$ such that $\nu^+:=\nu^+_\beta \in \mathcal{G}(\gamma^+)$, i.e. for all $\sigma \in \Omega$, for all $\Lambda \in \mathcal{S}$, 
$$
\nu^+[\sigma_\Lambda | \mathcal{F}_{\Lambda^c}](\omega')=\gamma_\Lambda^+(\sigma | \omega'),\; \nu^+{\rm -a.e.}(\omega').
$$
Now, proceeding like above in the proof of non-quasilocality, and using the definition of $\nu^+=T \mu^+$ as an image measure, one writes for $\nu^+$-a.e. $\omega'$,
$$
\nu^+[\sigma_\Lambda | \mathcal{F}_{\Lambda^c}](\omega')=\mu^+[\sigma_\lambda |\mathcal{F}_{\lambda^c}](\omega), \; {\rm with} \; \omega \in T^{-1}(\omega'), 
$$
where $\lambda^c:=(2\Z^2) \cap (2\Lambda)^c$, with $2\Lambda=\{j: \exists i \in \Lambda \; {\rm s.t.} \: j=2i\}$. When $\la \in \s$, $\lambda=(2\Z^2)^c \cup 2\Lambda$ is not finite so one has to use the global specification $\Gamma^+$ to write for $\nu^+$-a.e. $\omega'$:
$$
\nu^+[\sigma_\Lambda | \mathcal{F}_{\Lambda^c}](\omega')=\Gamma_\lambda^+(\sigma|\omega').
$$
Thus, one is led to introduce the following
\bd [Local specifications for the decimated measures]
Define for  $\Lambda \in \mathcal{S}$,
\be \label{gamma+-}
\gamma_\Lambda^{+}:=\Gamma^+_{(2\Z^2)^c \cup 2\Lambda} \; {\rm and} \;  
\gamma_\Lambda^{-}:=\Gamma^-_{(2\Z^2)^c \cup 2\Lambda}
\ee
\ed
so that we get
\begin{theorem}\cite{FP,FLNR1}\label{localspe}
The local specifications $\gamma^{+}$ and $\gamma^{-}$ composed by the kernels (\ref{gamma+-}) are monotonicity-preserving, $\gamma^{+}$ is right-continuous and $\gamma^{-}$ is left-continuous, but none of them is quasilocal at low enough temperature. The decimated measures of the extreme phases of the 2d-Ising model, $\nu^+=T \mu^+$ and $\nu^-=T \mu^-$, coincide with the weak limits
\be \label{select}
\nu^+(\cdot)=\lim_{\Lambda \uparrow Z^2} \gamma^{+}_\Lambda(\cdot \mid +) \; {\rm and} \; \nu^-(\cdot)=\lim_{\Lambda \uparrow \Z^2} \gamma^{-}_\Lambda(\cdot \mid -)
\ee
so that $\nu^{+} \in \mathcal{G}(\gamma^{+})$ and $\nu^{-} \in \mathcal{G}(\gamma^{-})$. They are extremal elements of the Choquet simplices $\mathcal{G}(\gamma^{+})$ and $\mathcal{G}(\gamma^{-})$ and for any $f$ increasing,  any $\nu \in \mathcal{G}(\gamma^+)$ or $\mathcal{G}(\gamma^-)$ $$\gamma^-_\Lambda f \leq \gamma_\Lambda^+ f \; {\rm  and} \; 
\nu^-[f] \leq \nu[f] \leq \nu^+[f].$$
\end{theorem}

Non-quasilocality is a direct consequence of the non-Gibbsian result of van Enter {\em et al.}, due itself to a possible hidden phase transition when conditioning w.r.t. the alternate configuration. This phenomenon has been extended in \cite{FP} to provide a general criterion.
\subsection{Fern\'andez-Pfister's criterion for quasilocality}

When a specification $\gamma^+$ is {\em monotonicity-preserving}, is is not difficult to see\footnote{See e.g. \cite{FP} p1304 or proceed like in p1294-1295.} that 
\be \label{conti}
\omega \in \Omega_{\gamma^+} \; \Longleftrightarrow \; \forall \Lambda \in \s, \lim_{\Delta \uparrow \Z^2} \gamma^+_\Lambda(\sigma | \omega_\Delta +_{\Delta^c})=\lim_{\Delta \uparrow \Z^2} \gamma^+_\Lambda(\sigma | \omega_\Delta -_{\Delta^c}) .
\ee

From this, and the expression (\ref{gamma+-}) of $\gamma^+$ in terms of the constraint measure (\ref{Constr}), Fern\'andez and Pfister deduce that uniqueness for the "constrained  specification" $\gamma^{(2\Z^2)^c,\omega}$ is a sufficient condition for quasilocality at $\omega$. Define the tail-measurable set
\be\label{omegpm}
\Omega_\pm:=\big\{ \omega : |\mathcal{G}(\gamma^{(2\Z^2)^c,\omega})|=1 \big\}.\ee
One first establishes that the latter is  the set of $\omega$'s where the specifications $\gamma^+$ and $\gamma^-$ coincide, and that this implies that it is contained in the sets of good configurations for $\gamma^+$ and $\gamma^-$, illustrating the paradigm "{\em Non-Gibbsianness is (mostly) due to an hidden phase transition}".
\begin{lemma}\cite{FP}\label{critql2}
The set of uniqueness for the constrained  specification is also
\be\label{omegapm}
\Omega_{\pm}=\big\{ \omega : \gamma_\Lambda^{+}(\sigma \mid \omega)=\gamma_\Lambda^{-}(\sigma \mid \omega), \forall \Lambda \in \mathcal{S} \big\}.
\ee
and it contains the set of continuity points $\Omega_{\gamma^+}$ and $\Omega_{\gamma^-}$: $\Omega_{\pm} \subset \Omega_{\gamma^{+}} \; {\rm and} \; \Omega_\pm \subset \Omega_{\gamma^{-}}$.
\end{lemma}
By monotonicity, one uses the absence of phase transition when the even sites are frozen in the $+$ and $-$-states\footnote{The freezing acts as an all positive/negative magnetic field, for which uniqueness holds from Theorem \ref{PTIsing}.} ($|\mathcal{G}(\gamma^{(2\Z^2)^c,+})|=|\mathcal{G}(\gamma^{(2\Z^2)^c,-})|=1 $) to fully identified the sets of continuity points of the decimated specifications with $\Omega_\pm$, and to eventually get the following criterion for the almost quasilocality of the decimated measures $\gamma^+$ and $\gamma^-$.
\begin{lemma}\cite{FP}\label{critql}{\bf (Criterion for Almost Gibbs)}
\begin{eqnarray} \label{crit}
\nu^{+}(\Omega_{\gamma^{+}})=1 \; &{\rm iff} &\; \nu^{+} \in \mathcal{G}(\gamma^{-})\\
 \nu^{-}(\Omega_{\gamma^{-}})=1 \; &{\rm iff}& \; \nu^{-} \in \mathcal{G}(\gamma^{+})\nonumber
 \end{eqnarray}
and $$\Omega_{\gamma^+}=\Omega_{\gamma^-}=\Omega_\pm = \big \{\omega : \gamma_\Lambda^+(\cdot|\omega)=\gamma_\Lambda^-(\cdot | \omega), \forall \Lambda \in \mathcal{S} \big\}.$$
\end{lemma}

This is a first step to prove that  $\nu^+$ and $\nu^-$ are {\em almost Gibbs}, but it is not enough. To check the rhs of the criterion (\ref{crit}), we  investigate in \cite{FLNR1} thermodynamic properties of the decimated measures and to establish a partial variational principle\footnote{A more complete variational principle has been established in \cite{KLNR}, but we do not need it here.}.

\subsection{Thermodynamic Properties and Variational Principles}

This section deals partially with the {\em second part of the Dobrushin program}, whose aim is to restore thermodynamics of renormalized Gibbs measures, in order to interpret them as equilibrium state\footnote{A complete interpretation and characterization in the original Gibbsian context can be also found in \cite{Is}.}. This part has been mostly achieved in \cite{FLNR1,KLNR} -- although previous partial statements exist in \cite{lef1} --  while the existence of thermodynamic functions (relative entropy, pressure) relies on results of Pfister \cite{P}. We do not focus on the latter thermodynamical approach but consider more statistical mechanical (specification-dependent) variational principles. They characterize equilibrium states - i.e. consistent measures for a given specification - in terms of the more local and physical property of {\em zero relative entropy}. As usual in the equilibrium approach to Gibbs measures, we focus in this part on {\em translation-invariant} objects. 

We introduce first the {\em
Kolmogorov-Sinai entropy}, well-defined for all $\mu \in \mathcal{M}_{1,\rm{inv}}^+$ as
\be \label{KS} h(\mu):=- \lim _{n \to \infty}
\frac{1}{|\la_n|}\sum_{\sigma_{\la_n}} \mu(\sigma_{\la_n})\log
\mu(\sigma_{\la_n}) \ee
and more generally, for $\mu, \nu \in \mathcal{M}_1^+$, the {\em relative entropy density} of $\mu$ w.r.t. $\nu$ to be the limit
\begin{equation}\label{eq:red}
h(\mu | \nu)\;=\;\lim_{n \to \infty} \frac{1}{|\la_n|}\sum_{\sigma_{\la_n} \in
\Omega_{\Lambda_n}} \mu(\sigma_{\la_n}) \cdot \log
\frac{\mu(\sigma_{\la_n})}{\nu(\sigma_{\la_n})}
\end{equation}
provided it exists. 
 It does for all 
$\mu\in \mathcal{M}_{1,\rm{inv}}^{+}$ when
$\nu\in\mathcal{M}_{1,\rm{inv}}^{+}$ is a Gibbs measure\footnote{Consistent with a t.i. UAC potential. It does not imply its existence for t.i. quasilocal measure, because the potential built by Kozlov from a t.i. quasilocal specification is not necessarily translation invariant \cite{KLNR}.} and, more generally, if $\nu$ is
asymptotically decoupled\footnote{Introduced
by Pfister \cite{P} to state general large deviation principles.}. This result is extended in  \cite{KLNR} for
general {\em translation-invariant quasilocal measures} and allows to characterize equilibrium states via the
\begin{definition}[Variational Principle relative to a Specification]\label{chp:defr}
 A variational principle occurs for a specification $\gamma$ and $\nu \in 
\mathcal{G}_{{\rm inv}}(\gamma)$
iff
\begin{equation}\label{eq:r4}
\forall \mu\in\mathcal{M}_{1,\rm{inv}}^{+}, \; h(\mu|\nu)
=0 \ \Longleftrightarrow\ \mu\in \mathcal{G}_{\rm inv}(\gamma)\;.
\end{equation}
\end{definition}
This result is well known for Gibbs measures consistent with a
translation-invariant UAC potential \cite{Ge} and has  been
extended recently to translation-invariant quasilocal DLR measures in \cite{KLNR}. More generally, in the right-continuous case, we have the  

\begin{theorem}\label{r.teo1}\cite{FLNR1}
Let $\gamma$ be a  {\em right-continuous} specification and $\nu\in\mathcal{G}_{\rm inv}(\gamma)$.
If $\mu\in\mathcal{M}_{1}^+$ is such that $h(\mu|\nu)=0$,
then
\begin{equation}\label{eq:r.10}
\mu\in \mathcal{G}(\gamma) \; \Longleftrightarrow \; \nu\Bigl[\,
g_{\Lambda_n\setminus\Lambda} \cdot \Bigl(\gamma_\Lambda f(\cdot_{\Lambda_n}+_{\Lambda_n^c})
- \gamma_\Lambda f(\cdot)\Bigr) \Bigr]
\;\mathop{\longrightarrow}\limits_{n\to\infty}\;0
\end{equation}
for all $\Lambda\in\mathcal{S}$ and $f\in\mathcal{F}_{\rm loc}$,
 where $g_{\Lambda_n \setminus \Lambda}:=\frac{d\mu_{\Lambda_n\setminus\Lambda}}{d\nu_{\Lambda_n\setminus\Lambda}}$
 provided it exists.
\end{theorem}

Thus, to get information about consistency from zero relative
entropy requires that the concentration properties of the density
of $\mu_{\Lambda_n \setminus \Lambda}$ w.r.t $\nu_{\Lambda_n
\setminus \Lambda}$ to beat asymptotic divergence due to the lack
of continuity of $\gamma$. When $\gamma$ is quasilocal, with $\Omega_\gamma=\Omega$, it  yields  the standard proof of this second part for translation-invariant quasilocal specifications.

\subsection{Generalized Gibbsianness for the 2d Decimated measure}
We describe now the two main notions proposed within the Dobrushin program, depending if one focuses on a relaxation of the
convergence properties of the potential ({\em Weak
Gibbsianness}), or
 on a relaxation of  topological properties ({\em Almost quasilocality} or {\em Almost Gibbsianness}). 
\subsubsection{Almost Quasilocality -- Almost Gibbsianness}
Recall that the set $\Omega_\gamma$ of {\em good configurations} of
a specification $\gamma$ gathers its points of continuity.
\begin{definition}[Almost Quasilocal (or Almost Gibbs)]
A probability measure $\mu \in \mathcal{M}_1^+$ is said to be {\em almost quasilocal} if there is a specification $\gamma$ s.t. $\mu \in \mathcal{G}(\gamma)$ and $\mu(\Omega_\gamma)=1$.
\end{definition}

We start from the specification $\gamma^{I}=\gamma^{\beta\Phi}$ defined on $\mathbb{Z}^2$ by (\ref{BG}) for the potential (\ref{ising}), at low enough temperature $\beta^{-1}$ in such a way that the decimation $\nu^+=\nu_\beta^+$ of the $+$-phase $\mu^+_\beta$ is not quasilocal. Then, putting together Lemma \ref{critql} and Theorem \ref{r.teo1}, we get

\begin{theorem}\label{AGChoquet}[{\bf Almost quasilocality of the decimated measures}] Let $\nu^+=T \mu^+$ and $\nu^-=T \mu^-$ be the decimation of the extremal phases of the 2d-Ising model at {\em any} temperature.
\ben
\item $\nu^+$ (resp. $\nu^-$) is {\em almost quasilocal}: it is consistent with the specification $\gamma^+$ (resp. $\gamma^-$) whose set of continuity points is of full measure : $\nu^+(\Omega_{\gamma^+})=1$ (resp. $\nu^-(\Omega_{\gamma^-})=1$).

\item In fact, the set $\Omega_\pm$ (\ref{omegapm}), where the specifications $\gamma^-$ and $\gamma^+$ coincide, is of full ($\nu^-$ and $\nu^+$) measure and coincides with the sets of continuity points $\Omega_{\gamma^+}$ and $\Omega_{\gamma^-}$. Moreover
\be \label{choquet}
\mathcal{G}(\gamma^+)=\mathcal{G}(\gamma^-)=[\nu^-,\nu^+]
\ee
is a Choquet simplex whose extremal elements are $\nu^-$ and $\nu^+$, such that for any $f$ bounded monotone increasing, any $\nu \in \mathcal{G}(\gamma^+)$, we have $\nu^-[f] \leq \nu[f] \leq \nu^+[f].$
\item These measures are also equilibrium states in the sense that, for e.g. the $+$-phase:
\ben
\item For all $\nu \in \mathcal{M}_{1,{\rm inv}}^+$, the relative entropy $h(\nu | \nu^+)$ exists. In particular $h(\nu^- | \nu^+)=0$.
\item If $h(\nu | \nu^+)=0$ and $\nu(\Omega_\pm)=1$ for $\nu \in \mathcal{M}_{1}^+$, then $\nu \in \mathcal{G}(\gamma^+)$ is also almost Gibbs.
\een
\een
\end{theorem}
{\bf Proof of Theorem \ref{AGChoquet}}
\ben
\item
 By Lemma \ref{critql}, almost quasilocality of $\nu^+$ is equivalent of having $\nu^- \in \mathcal{G}(\gamma^+)$, which is proved now using Theorem \ref{r.teo1}. First, a standard result\footnote{See e.g. \cite{VEFS} p 970 for a one-line proof.} yields 
$$
h(\nu^- \mid \nu^+)=h(T\mu^- \mid T\mu^+) \leq h(\mu^- \mid \mu^+)=0.
$$
The lhs of (\ref{eq:r.10}) is true by right-continuity of $\gamma^+$, so that one indeed gets $\nu^- \in \mathcal{G}(\gamma^+)$. 

\item  We know from Lemma \ref{critql} that $\Omega_\pm=\Omega_{\gamma^+}=\Omega_{\gamma^-}$ is of full $\nu^+/\nu^-$-measure. Using the definition of  $\gamma^+$  in terms of the global specification $\Gamma^+$ and that  $\mu^-$ and $\mu^+$ are weak limits obtained with all $-$ or all $+$-boundary conditions, one gets
$$
\nu^-(\cdot)=\lim_{\Lambda \uparrow \mathcal{S}} \gamma_\Lambda^+(\cdot \mid -) \; \rm{and} \; \nu^+(\cdot)=\lim_{\Lambda \uparrow \mathcal{S}} \gamma_\Lambda^+(\cdot \mid +).
$$
For $f$ increasing and $\nu \in \mathcal{G}(\gamma^+)$, $\limsup_\Lambda \gamma_\Lambda^+ [f | \cdot]$ is a version of $\nu [f | \mathcal{F}_\infty](\cdot)$. By monotone-preservation, one has, $\nu$-a.s., $\nu^-[f] \leq \nu [f |\mathcal{F}_\infty](\cdot) \leq \nu^+[f]$
and thus
\be \label{fkg}
\nu^-[f] \leq \nu[f] \leq \nu^+[f].
\ee
The fact that $\mathcal{G}(\gamma^+)$  is a Choquet simplex is a general property of any specification \cite{Dyn, ALN2}, and (\ref{fkg}) implies that $\nu^-$ and $\nu^+$ are its extremal elements, so that $\mathcal{G}(\gamma^+)=[\nu^-,\nu^+]$. Proceeding similarly for  $\gamma^-$, one also gets 
$$
\nu^-(\cdot)=\lim_{\Lambda \in \mathcal{S}} \gamma_\Lambda^-(\cdot \mid -) \; \rm{and} \; \nu^+(\cdot)=\lim_\Lambda \gamma_\Lambda^-(\cdot \mid +)
$$
and (\ref{fkg}) for all $\nu \in \mathcal{G}(\gamma^-)$. So we get $\mathcal{G}(\gamma^-)=[\nu^-,\nu^+]=\mathcal{G}(\gamma^+)$.
\item
It is proved in \cite{FLNR1} as a consequence of Lemma \ref{r.teo1}.

\een
\subsubsection{Weak Gibbsianness}

\begin{theorem}[Weak Gibbsianness of the decimated measures]\cite{MRSM}\label{thm9}
The decimated measure $\nu^+$ is weakly Gibbs at any temperature: there exists an a.s. absolutely convergent potential consistent with the specification $\gamma^+$. Moreover, there exists a t.i. potential $\Psi=\Psi^+$  absolutely convergent on a set $\Omega_{\Psi^+}$ of full $\nu^+$-measure,  s.t. $\nu^+$ is consistent with $\Psi$.
\end{theorem}

{\bf Proof of Theorem \ref{thm9}:} The first tool  is to use the right-continuity  to construct first a convergent {\em vacuum potential} consistent with the specification $\gamma^+$, written $\gamma$ for simplicity.\\

{\bf Step 1: Convergence and consistency of the vacuum potential}\\

Let us  assume for
the moment that a vacuum  potential $\Phi^+$, with vacuum
state + exists. Then, for any $\Lambda \in
\mathcal{S}$, one has  $H_\Lambda^{\Phi^+}(+ |+)=0$, which implies $\gamma_\Lambda(+|+)=Z_\Lambda (+)^{-1} > 0$, 
so for  any  $\sigma \in \Omega$, $\gamma_\Lambda(\sigma|+)=\gamma_\Lambda(+|+) \; e^{-
H_\Lambda^{\Phi^+}(\sigma | +)}.$
Then, by non-nullness of the specification,
\be \label{relat}
H_\Lambda^{\Phi^+}(\sigma | +)=-\ln
\frac{\gamma_\Lambda(\sigma|+)}{\gamma_\Lambda(+|+)}.
\ee
\begin{equation}\label{Hambc0}
H_\Lambda^{\Phi^+}(\sigma | +)=
 \sum_{A \cap \Lambda \neq \emptyset} \Phi^+_A(\sigma_\Lambda
+_{\Lambda^c})=\sum_{A \subset \Lambda} \Phi^+_A (\sigma) +
\sum_{A \cap \Lambda \neq \emptyset,A \cap \Lambda^c \neq
\emptyset} \Phi^+_A(\sigma_\Lambda +_{\Lambda^c}).
\end{equation}
The last sum is null by the vacuum property: prescribing the vacuum state as a boundary condition is equivalent to consider free b.c. and (\ref{relat}) yields
\begin{equation}\label{Hfree2}
\forall \Lambda \in \mathcal{S},\; \forall \sigma \in \Omega,\;
H_\Lambda^{\Phi^+,f}(\sigma)=
\sum_{A \subset \Lambda} \Phi^+_A(\sigma)=-\ln
\frac{\gamma_\Lambda(\sigma|+)}{\gamma_\Lambda(+|+)}.
\end{equation}
To recover the potential $\Phi_A^+$ for any set A, one proceeds now by induction (as in  \cite{Fer,ALN2}) using 
\begin{proposition}[Moebius  inversion
formula] Let $\mathcal{S}$ be a countable set of finite sets,
$H=(H_\Lambda)_{\Lambda \in \mathcal{S}}$ and $\Phi=(\Phi_A)_{A
\in \mathcal{S}}$ be set functions from $\mathcal{S}$ to
$\mathbb{R}$. Then
\begin{equation}\label{Moeb}
\forall \Lambda \in \mathcal{S}, H_\Lambda=\sum_{A \subset
\Lambda} \Phi_A \; \Longleftrightarrow \; \forall A \in
\mathcal{S}, \; \Phi_A= \sum_{B \subset A} (-1)^{|A \setminus B |}
H_B.
\end{equation}
\end{proposition}
\begin{proposition}
The potential $\Phi^+$ defined for all $\sigma \in \Omega$ and for all $A \in \mathcal{S}$ by
\begin{equation}\label{vacuumpot}
\Phi^+_A(\sigma) = - \sum_{B \subset A} (-1)^{|A \setminus B |}
\ln \frac{\gamma_B(\sigma | +)}{\gamma_B(+ | +)}.
\end{equation}
is a vacuum potential with vacuum state $+$. It is moreover convergent and consistent with $\nu^+$.
\end{proposition}
It is obviously a potential and we prove first that $\Phi^+$ satisfies  the vacuum
condition. Write
\begin{equation}\label{vacpot}
\Phi^+_A(\sigma)=-\sum_{B \subset A} (-1)^{|A \setminus B|}\ln
\frac{\gamma_B(\sigma | +)}{\gamma_B(+ | +)}=\sum_{B \subset A}
(-1)^{|A \setminus B|} H_B^{\Phi^+,f} (\sigma)
\end{equation}
where by the Moebius formula (\ref{Moeb}), one has for all $B \in \mathcal{S}$, $H_B^{\Phi^+,f} (\sigma)=\sum_{A
\subset B} \Phi^+_A(\sigma)$.
Consider $A \in
\mathcal{S}$ and  $\sigma \in \Omega$ such that there exists $i
\in A$ where $\sigma_i=+_i$. Using consistency of the specification $\gamma$ via the key-bar displacement property (\ref{Keybar}), one gets
$$
\frac{\gamma_B(\sigma |+)}{\gamma_B(+ | +)}=\frac{\gamma_{B \setminus i}(\sigma |+)}{\gamma_{B \setminus i} (+ | +)}
$$
so that
\begin{equation}\label{consistenciii}
\forall i \in B \subset A, \; H_B^{\Phi^+,f} (\sigma)= H_{B
\setminus i}^{\Phi^+,f} (\sigma)
\end{equation}
which implies the vacuum property -- crucial to get consistency -- using expression (\ref{vacpot}).\\

To check that $\Phi^+$ is convergent and consistent with $\gamma$, we first need to be able to 
define
$$
\forall \sigma,\omega \in  \Omega, \; H_\Lambda^{\Phi^+} (\sigma |
\omega)=\sum_{A \cap \Lambda \neq \emptyset, A \in \mathcal{S}}
\Phi^+_A(\sigma_\Lambda \omega_{\Lambda^c})
$$
to extend  the definition of the Hamiltonian with free
b.c. (\ref{Hfree2}) to an Hamiltonian with b.c.  $\omega \in \Omega$. It amounts to prove the convergence of
the potential, in the sense that for all $\sigma \in \Omega$,
 the limit as $\Delta \uparrow \mathcal{S}$ of
the net $\Big(\sum_{A \cap \Lambda \neq \emptyset, A \subset
\Delta} \Phi^+_A(\sigma) \Big)_{\Delta \in \mathcal{S}}$ is
finite. To do so, write 
\begin{displaymath}
\sum_{A \cap \Lambda \neq \emptyset,A \subset
  \Delta}\Phi_A^+(\sigma)=\sum_{A \subset \Delta}\Phi_A^+(\sigma)-\sum_{A
  \subset \Delta \cap \Lambda^c}\Phi_A^+(\sigma)=-\ln{\frac{\gamma_\Delta(\sigma|+)}{\gamma_{\Delta}(+|+)}}+
 \ln{\frac{\gamma_{\Delta \cap \Lambda^c}(\sigma|+)}{\gamma_{\Delta \cap \Lambda^c}(+|+)}}.
\end{displaymath}
To deal with the second term, we use first the properties of the specification to write that $\gamma_{\Delta \cap \Lambda^c} (\sigma |+)=\gamma_{\Delta \cap \Lambda^c} (+_\Lambda \sigma_{\Lambda^c}|+)$, in order to use consistency (\ref{Keybar}) on $\Delta \cap \Lambda^c$ and $\Delta$, to get
$$
\frac{\gamma_{\Delta \cap \Lambda^c}(\sigma|+)}{\gamma_{\Delta \cap \Lambda^c}(+|+)}=\frac{\gamma_{\Delta \cap \Lambda^c}(\sigma_{\Delta \cap \Lambda^c} +_{\Delta^c \cup \Lambda} |+)}{\gamma_{\Delta \cap \Lambda^c}(+|+)}=\frac{\gamma_{\Delta}(\sigma_{\Delta \cap \Lambda^c} +_{\Delta^c \cup \Lambda} |+)}{\gamma_{\Delta }(+|+)}
$$
so that
\be \label{equin}
\sum_{A \cap \Lambda \neq \emptyset,A \subset
  \Delta}\Phi_A^+(\sigma)=-\ln{\frac{\gamma_{\Delta}(\sigma_{\Delta
  \cap \Lambda^c}+_{\Delta^c \cup \Lambda}| +)}
{\gamma_{\Delta}(\sigma | +)}}=-\ln \frac{\gamma_{\Lambda \cap \Delta}(\sigma_{\Delta}+_{\Delta^c}|\sigma_{\Delta}+_{\Delta^c})}{\gamma_{\Lambda \cap \Delta}(\sigma_{\Delta \cap  \Lambda^c} +_{\Delta^c \cup \Lambda}|\sigma_{\Delta \cap  \Lambda^c} +_{\Delta^c \cup \Lambda})}
\ee
using again (\ref{Keybar}) for
$\sigma_{\Delta
  \cap \Lambda^c}+_{\Delta^c \cup \Lambda}$ and $\sigma_{\Delta}+_{\Delta^c}$ which agree on $\Delta \cap
  \Lambda$.
Let $\Delta \uparrow \mathcal{S}$ in the sense defined, big enough so that $\Delta \supset
\Lambda$ and thus for any  $\sigma \in \Omega$
\be \label{partdelt}
\sum_{A \cap \Lambda \neq \emptyset,A \subset
\Delta}\Phi_A^+(\sigma)=-\ln \frac{\gamma_{\Lambda}(\sigma_{\Delta}+_{\Delta^c}|\sigma_{\Delta}+_{\Delta^c})}{\gamma_{\Lambda}(\sigma_{\Delta \cap  \Lambda^c} +_{\Delta^c \cup \Lambda}|\sigma_{\Delta \cap  \Lambda^c} +_{\Delta^c \cup \Lambda})}=-\ln \frac{\gamma_\Lambda(\sigma | \sigma_\Delta +_{\Delta^c})}{\gamma_\Lambda(+|\sigma_\Delta +_{\Delta^c})}.
\ee
Now, {\bf by right-continuity}, we get convergence of $\Phi^+$ and that for all $\sigma \in \Omega$,
\be \label{Ham}
H_\Lambda^{\Phi^+}(\sigma)=\sum_{A \cap \Lambda \neq
  \emptyset, A \in \mathcal{S}}\Phi_A^+(\sigma)=-\ln \frac{\gamma_\Lambda(\sigma | \sigma)}{\gamma_\Lambda(+ | \sigma)}< + \infty.
\ee
Using the expression (\ref{Hambc0}) of the Hamiltonian with boundary condition $\omega$, we also get
\be \label{Hambc2} \forall \omega \in
\Omega, \; H_\Lambda^{\Phi^+}(\sigma | \omega)=- \ln
\frac{\gamma_\Lambda (\sigma | \omega)}{\gamma_\Lambda (+ |
\omega)}< + \infty \ee
so that $\gamma_\Lambda(\sigma | \omega) = \gamma_\Lambda(+ | \omega) \; e^{-H_\Lambda^{\Phi^+}(\sigma | \omega)}$
and $\Phi^+$ is also consistent with $\gamma$.\\

{\bf Step 2 : Construction of a consistent telescoping  potential}\\

In fact,  the vacuum potential with vacuum state $+$ defined above is convergent and consistent with $\gamma$ {\em if and only if} the latter is right-continuous \cite{MRVM}.
Nevertheless, even when $\gamma$ is quasilocal, this vacuum potential is not UAC and is thus not enough to recover a proper Gibbsian property
(\ref{UAC}). To gain summability and absoluteness,
 Kozlov \cite{Ko} introduced a particular re-summation procedure by telescoping the terms of the Hamiltonian with
 free boundary conditions in annuli large enough to recover absoluteness, carefully keeping consistency,
  to eventually get a potential $\Psi$ such that (\ref{UAC}) holds and, for all $\sigma \in \Omega$
\begin{equation}\label{conskoz}
\sum_{A \subset \Lambda}
\Psi_A(\sigma)=H_\Lambda^{\Phi^+,f}(\sigma)=\sum_{A \subset
\Lambda} \Phi^+_A(\sigma) = -\ln \frac{\gamma_\Lambda(\sigma | +)}{\gamma_\Lambda(+ | + )}.
\end{equation}
Inspired from this procedure, also described in \cite{Fer,ALN2}, Maes {\em et al.} \cite{MRSM,MRVM,MRVM2}  introduced a general telescoping scheme and conditions to get an a.s. {\em absolutely} convergent and consistent potential. The main tool  is to introduce a family $(L_{i,m})_{i \in \Z^2, m \geq 0}$ s.t.
\be \label{onetoone}
\forall A \in \s, \exists ! (i,m) \; {\rm s.t.} A \ni i, A \subset L_{i,m}, A \cap L_{i,m-1}^c \neq \emptyset \; {\rm for} \; m \geq 1.
\ee
For any couple $(i,m) \in \Z^2 \times \N$, denote by $\mathbb{L}_{i,m}$  the set of $A$'s that satisfy (\ref{onetoone}). It is then possible  for all $\Lambda \in \s$ 
to rewrite formally\footnote{To be true, this re-writing requires absolute convergence, checked next section.}  the Hamiltonian (\ref{Ham}) as
\be \label{Reorder}
H_\Lambda^{\Phi^+}=\sum_{A \cap \Lambda \neq
  \emptyset, A \in \mathcal{S}}\Phi_A^+=\sum_{i \in \Z^2} \sum_{m \geq 0} \sum_{A \in \mathbb{L}_{i,m}, A \cap \Lambda \neq \emptyset} \Phi^+_A
\ee
and to define  the potential $\Psi$  from the vacuum potential by telescoping it in the annulus $L_{i,m}$, i.e. by putting $\Psi_A(\omega)=0$ unless $A=L_{i,m}$ with
\be \label{PotPsi0}
\Psi_{L_{i,m}}=\sum_{R \in \mathbb{L}_{i,m}} \Phi^+_R.
\ee
Using consistency (\ref{Keybar}) and (\ref{partdelt}) like in (\ref{equin}),  (\ref{PotPsi0}) becomes, for all $i \in \Z^2$ and $m \geq 1$,
\begin{eqnarray} \label{PsiLnk}
 \nonumber \Psi_{L_{i,m}}(\omega)&=& \sum_{R \ni i, R \subset L_{i,m}} \Phi^+_R(\omega) - \sum_{R \ni i, R \subset L_{i,m-1}} \Phi^+_R(\omega)\\
&=& - \ln \frac{\gamma_{\{i\}}(\omega | \omega_{L_{i,m}} +_{L_{i,m}^c})}{\gamma_{\{i\}}(+ | \omega_{L_{i,m}} +_{L_{i,m}^c})} + \ln \frac{\gamma_{\{i\}}( +| \omega_{L_{i,m-1}} +_{L_{i,m-1}^c})}{\gamma_{\{i\}}(\omega | \omega_{L_{i,m-1}} +_{L_{i,m-1}^c})}
\end{eqnarray}
while for $m=0$ one already gets (\ref{defpsiLim}). Using again (\ref{Keybar}) twice and the notation $\omega^L:=\omega_L +_{L^c}$, (\ref{PsiLnk}) can be formally rewritten, for all $m \geq 1$,
\be \label{PotPsi3}
\Psi_{L_{i,m}}(\omega)=-\ln \frac{\gamma_{L_{i,m}}(\omega^{L_{i,m}} | +)\cdot \gamma_{L_{i,m}}(\omega^{L_{i,m-1}\setminus i}  | +)}{\gamma_{L_{i,m}}(\omega^{L_{i,m} \setminus i} | +) \cdot  \gamma_{L_{i,m}}(\omega^{L_{i,m-1}} |+)}.\\
\ee

\newpage

{\bf Step 3 : Absolute convergence and consistency of the telescoping potential}

It is a consequence of almost Gibbsianness, and more precisely of almost-quasilocality in any direction\footnote{See e.g. Proposition 4.24 in \cite{Fer}. Here, we have even more because the specification is right-continuous.}. Define for any $i \in \Z^2$ and $n \in \N$,
$$
g_i(n,\omega)=|\gamma_i(\omega |\omega) - \gamma_i (\omega |\omega_{\la_n} +_{\la_n^c}) |
$$
and $g(n,\omega):=\sup_i g_i(n,\omega)$.
The set $\Omega_\gamma^+$ of points of right-continuity 
$$
\Omega_\gamma^+=\big\{\omega \in \Omega : g(n,\omega):\longrightarrow_n 0 \big\}
$$
satisfies thus $\nu^+(\Omega_\gamma^+)=1$ (in fact here $\Omega_\gamma^+=\Omega$). To use non-nullness, define also $m_i:=\inf_{\omega \in
\Omega}\gamma_{\{i\}} (\omega | \omega) >0$ and use the standard inequality
\be \label{standineq}
\forall a,b >0,\; | \ln(a) - \ln(b) | \leq \frac{| a - b |}{\min(a,b)}
\ee
to get from (\ref{partdelt}) the bound
\be \label{bound}
\forall n \in \N, \forall i \in \Z^2, \forall \omega \in \Omega_\gamma^+ ,\; \big| \sum_{B \ni i, B \subset \la_n} \Phi^+_B(\omega) \big| \leq \frac{g(n,\omega)}{m_i}.
\ee
Now, order lexicographically $\Z^2$ as $\{x_1,x_2,\dots \}$ and consider  $(n_l)_{l \in \N}$ s.t., for all $\omega \in \Omega_\gamma^+$,
\be \label{subseq}
\sum_{l=0}^{\infty} g(n_l,\omega) < \infty,
\ee
and eventually choose the sets $L_{i,m}$ to coincide with rectangles around each of the site $x_l$ and to satisfy the condition (\ref{onetoone}) with $m=n_l$. From the definition (\ref{PotPsi0}) of the potential, we get the absolute convergence on $\Omega_\gamma^+$ thanks to the bound
$\big| \Psi_{L_i,m}(\omega) \big| \leq \frac{g(n_l,\omega) + g(n_{l-1},\omega)}{m_{n_l}}.$ This absolute convergence justifies the re-summation procedures and consistency is kept.

\section{Parsimonious Description of the Decimated Measures}

We consider now our decimated measures $\nu^+$ and $\nu^-$ -- together with their specifications $\gamma^+$ and $\gamma^-$ -- under the light of  Parsimony and study whether they could provide  another non-Gibbsian example of such fields after this of \cite{CGL}.\\

Recall first the expression got for the specification $\gamma^+$ in terms of the constrained measures $\mu^{+,\omega}_{S}$  got themselves as weak limits of the constrained specifications $\gamma^{S,\omega}$ with $+$-boundary conditions. By (\ref{gamma+-}), one has indeed, for all $\Lambda \in \mathcal{S}$, for all $\omega \in \Omega$,
$$
\gamma_\Lambda^+(\cdot \mid \omega)=\Gamma^+_S (\cdot \mid \omega)=
\lim_{\Lambda \uparrow S} \gamma_\Lambda^I(\cdot \mid +_S \omega_{S^c})
$$

where $\gamma^I$ is the original Ising specification and the subset of non-even sites is $S=(2 \mathbb{Z}^2)^c \cup 2 \Lambda$. Another way to understand this limit is to consider first the weak limit (\ref{Constr}) so that
\begin{equation}\label{ParsSpe}
\forall \Lambda \in \s,\; \gamma_\Lambda^+(d\eta \mid \omega) = \mu_S^{+, \omega} (d\eta_S) \otimes \delta_{\omega_{S^c}}(d \eta_{S^c}).
\end{equation} 

To identify the family of contexts $(C_\Lambda(\omega))_{\Lambda,\omega}$ as in Definitions {\ref{DefParsRand} and {\ref{DefParsMark}, we would need to investigate the dependency of (\ref{ParsSpe}) in $\omega$, for a given $\Lambda$,  and thus to investigate the dependency in $\omega$ of the support of the constrained measure $\mu_S^{+,\omega}$. Even at high temperature  it seems that Definition {\ref{DefParsRand} is still too strong to incorporate decimated measures in this framework, while at low temperature one could expect even more complicated dependencies. Indeed, in this uniqueness high temperature region, for any boundary condition $\omega \in \Omega$, one has $\mathcal{G}(\gamma^{S,\omega}) = \big\{ \mu_{S^c}^\omega \big\}$ and for e.g. the $\omega=+$ b.c., one gets for $\gamma_\Lambda^+(\cdot|+)$ an Ising model on the decorated lattice with magnetic fields acting  on some even sites only, with strength 2. The constrained measure being an infinite volume Gibbs measures, it might depend on everything outside $\Lambda$ and in particular one should have  $\mid C_\Lambda(\omega) \mid =\infty$ for many $\omega's$ (for example, $C_\Lambda(+)$ seems to be of the order of $\Lambda^c$). We do not get any explicit expression for the contexts, nor a characterization of their typical size and do not formally prove here that decimated measures are not Parsimonious Markov Fields, nor Parsimonious Random fields, but we show that an asymptotic Parsimonious description is still possible if one adapts the definition and allow weak dependencies for more distant regions.

We enlarge thus this parsimonious context and describe our decimated measures in a weakly Gibbsian form with an a.s. absolutely convergent potential that admits a quenched exponential decay of correlation. The hard part of the job has been done in a more general context by Maes {\em et al.} \cite{MRSM,shl}. Such a configuration-dependent description requires to focus on typical properties of the phase to characterize the structure of configurations that could create long dependencies in the conditioning. Path large deviation (PLD (\ref{PLD})) and Amoeba (Definition \ref{AmoebaProp}), designed for that purpose, are such properties. The fact that the latter is typical means that $\mu_\beta^+$-a.s., beyond some configuration-dependent distance,  the density of minuses is small enough to make the far away contributions exponentially small. The configuration-dependent distance is then interpreted as a Quenched correlation decay (QCD(\ref{Corrlenght}) below) and will be our tool to get a parsimonious description of the transformed measure. 

\begin{theorem}[Parsimonious description via weak Gibbs]\label{thm10}
Consider the weakly Gibbs\-ian description of Theorem \ref{thm9} for the decimated measure $\nu^+$. It is possible to get it with an absolutely convergent potential $\Psi$ s.t. t $\Psi_A$ is non-null only for a set $A$ of the form 
$$
L_{i,m}= \big\{k \in \Z^2 : k\leq i, \; ||k-i||_1 \leq m \big \}, \; {\rm for \; any} \; i \in   \Z^2,\; m \in \N.
$$
It is defined in terms of (\ref{gamma+-}), the specification $\gamma=\gamma^+$ : for all $i \in \Z^2$ and $m \geq 1$,
\begin{eqnarray} \label{defpsiLim}
\Psi_{L_{i,0}}(\omega)&=&-\ln \frac{\gamma_{ \{i \}}(\omega | \omega_{L_{i,0}} +_{L_{i,0}^c})}{\gamma_{ \{i \}}(+ | \omega_{L_{i,0}} +_{L_{i,0}^c})}\\
\Psi_{L_{i,m}}(\omega)&=&-\ln \frac{\gamma_{L_{i,m}}(\omega^{L_{i,m}} | +)\cdot \gamma_{L_{i,m}}(\omega^{L_{i,m-1}\setminus i}  | +)}{\gamma_{L_{i,m}}(\omega^{L_{i,m} \setminus i} | +) \cdot  \gamma_{L_{i,m}}(\omega^{L_{i,m-1}} |+)}
\end{eqnarray}
Moreover, at low enough temperature, for all $\omega \in \Omega_{\Psi^+}$ and $i \in \Z^2$, there exist {\em configuration-dependent lengths} $l_i(\omega)$, positive constants $C_1,C_2 < + \infty$ and $0<\lambda <\frac{1}{2}$ s.t., when $A=L_{i,m}$,
\be \label{Corrlenght}
|\Psi^+_A(\omega)| \leq C_1 \cdot \mathbf{1}_{m \leq l_i(\omega)} \cdot m + C_2 \cdot \mathbf{1}_{m > l_i(\omega)} \cdot m \cdot e^{-\lambda m}.
\ee
Similar results hold for  $\nu^-$, with a different a.s. absolutely convergent potential $\Psi^-$.
\end{theorem}

The "quenched"  lengths $l_i(\omega)$, given in the course of the 
proof, correspond to the distance beyond which the amoebas compatible with 
$\xi=\xi(\omega)$ are benign, according to Definition \ref{AmoebaProp}. The 
results of \cite{shl} (Theorem \ref{Amoebathm}) yield that these lengths are typically finite and the 
configurations where it holds will be the points of convergence of the 
potential. The latter is constructed by a telescoping procedure inspired by 
Kozlov's techniques for Gibbs measures and generalized by Maes {\em et al.} 
\cite{MRSM,shl}. It is expressed in terms of  correlations functions of a 
constrained measure $\mu^\xi$, got from $\mu^+$ by constraining the even 
sites to be under the configuration $\xi \in \Omega_2$ -- chosen  $\mu_\beta^+$-typical so that the 
amoeba property holds. This property is indeed enough to get a "quenched" 
contour estimate for the constrained measure $\mu^\xi$, which allows to get a 
QCD by an adaptation of the percolation techniques of Burton and Steiff \cite
{BS}. The choice of the sets $L_{i,m}$ allows  to tune the convergence 
properties of the potential for typical configurations. Kozlov \cite{Ko} also 
played on it to introduce two different potentials: an  absolutely convergent 
potential when the specification is quasilocal, and a {\em  t.i.} absolutely 
convergent potential under a condition stronger than quasilocality. In \cite
{MRVM}, this last condition is required for typical configurations, in order to 
get an {\em a.s. t.i.} absolutely convergent potential, while at low 
temperature we  get it as a consequence of the {\em Quenched Correlation 
Decay} (QCD) (\ref{Corrlenght}), following \cite{MRSM}. More precisely, 
consider for all $i \in \Z^2$ and $m \geq 0$, the sets
\be \label{Lim}
 L_{i,m}=\{k \in  \mathbb{Z}^2: k\leq i, ||j - i||_1 \leq m\} 
\ee
which satisfies the condition (\ref{onetoone}) and yields a t.i.\footnote{This is not the case if one removes the condition $k \leq i$.} potential $\Psi$ defined by (\ref{PsiLnk}) and (\ref{PotPsi3}).

To proceed and control its convergence, we shall work  at finite volume $\la_n$  first and denote by $\nu_n^+(\omega)$ the $\nu^+$-probability that the configuration coincide with $\omega$ on $\Lambda_n$, and is arbitrary outside. Use (\ref{select}) to get an approximation of $\Psi$ by a potential $\Psi^n$ expressed in terms of $\nu_n^+$. The behavior of $\Psi_n$ as $n$ grows will be controlled using correlations of the constrained measures and this will define the potential $\Psi$ properly, coinciding with (\ref{PotPsi3}). We introduce also $\xi^j \in \Omega_2$ with $(\xi^j)_l=\xi_l$ if $l \neq j \in 2\Z^2$ and $(\xi^j)_j=+1$ and write
$$
\Psi^n_{L_{i,0}}(\omega)=- \ln \frac{\nu_n^+(\omega^i)}{\nu_n^+(+)}=- \ln \frac{\mu_n^+(\omega^{2i})}{\mu_n^+(+)}=-\ln \frac{\mu_n^{\xi^j}(\omega)}{\mu_n^{\xi^j}(+)}
$$
and
 \be \label{PotPsi4}
\Psi_{L_{i,m}}^n (\omega) = - \ln \frac{\nu_n^+(\omega^{L_{i,m}})}{\nu_n^+(\omega^{L_{i,m-1}})} \cdot \frac{\nu^+_n(\omega^{L_{i,m-1} \setminus i})}{ \nu_n^+(\omega^{L_{i,m} \setminus i})}= - \ln \frac{\mu_n^+(\omega^{2L_{i,m}})}{\mu_n^+(\omega^{2L_{i,m-1}})}  \frac{\mu_n^+(\omega^{2L_{i,m-1} \setminus 2i})}{ \mu_n^+(\omega^{2L_{i,m} \setminus 2i})}.
\ee
To express it in terms of  the constrained measure, we further telescope the potential in $L_{i,m}$ for fixed $i \in \Z^2$ and $m \geq 1$. We write $v(i,m)$ for the size of the annulus $L_{i,m} \setminus L_{i,m-1}$, bounded below by $2m$. Order lexicographically this annulus, write it $\{ j_1,j_2,\dots,j_r,\dots, j_v\}$ and eventually introduce for $r=1,\dots,v$ the sets
$$
Q_{i,m,0}=L_{i,m-1}, \; Q_{i,m,v}=L_{i,m} \; \; {\rm and} \; \; Q_{i,m,r}:= L_{i,m-1} \cup \{ j_1,j_2,\dots,j_r\}.
$$
Provided $m>0$ and $v(i,m) \geq 1$, we can now telescope (\ref{PotPsi4}) as
\be \label{covariance0}
\Psi_{L_{i,m}}^n (\omega) = - \sum_{r=1}^{v(i,m)}  \ln \frac{\mu_n^+(\omega^{2Q_{i,m,r}}) \mu_n^+(\omega^{2Q_{i,m,r-1} \setminus 2i})}{\mu_n^+(\omega^{2Q_{i,m,r-1}}) \mu_n^+(\omega^{2Q_{i,m,r} \setminus 2i})}=- \sum_{r=1}^{v(i,m)}  \ln\frac{ \frac{\mu_n^+(\omega^{2Q_{i,m,r-1} \setminus 2i})}{ \mu_n^+(\omega^{2Q_{i,m,r}})}}{\frac{\mu_n^+(\omega^{2Q_{i,m,r-1}})}{\mu_n^+(\omega^{2Q_{i,m,r}})}  \frac{\mu_n^+(\omega^{2Q_{i,m,r} \setminus 2i})}{\mu_n^+(\omega^{2Q_{i,m,r}})}}.
\ee

Introduce the constrained measures $\mu_{n,+}^{\omega^{2Q}}$ with freezing in $\omega$ in $2 Q$, with "$+$"-boundary condition outside $\Lambda_n$. Using the Gibbs (n.n.) property, one gets for all $i$ and $k$ in $\Z^2$
$$
\mu_n^+(\omega^{2Q \setminus \{2k\}})=\mu_n^+(\omega^{2Q}) \cdot \mu_{n,+}^{\omega^{2Q}}[e^{2 \beta \sigma_{2i}}].
$$
For two close sites s.t. $|i-j_r| \leq 4$, we also get
$$
\mu_n^+(\omega^{2Q \setminus \{2i, 2j_r\}})=\mu_n^+(\omega^{2Q}) \cdot \mu_{n,+}^{\omega^{2Q}}[e^{2 \beta (\sigma_{2i}\sigma_{2j_r})}]
$$
or in a factorized form when $|i-j_r|>4$
$$
\mu_n^+(\omega^{2Q \setminus \{2i, 2j_r\}})=\mu_n^+(\omega^{2Q}) \cdot \mu_{n,+}^{\omega^{2Q}}[e^{2 \beta \sigma_{2i}}e^{2 \beta \sigma_{2j_r}}].
$$
Taking $Q=Q_{i,m,r}$, one can rewrite 
\be \label{covariance1}
\Psi_{L_{i,m}}^n (\omega) =- \sum_{r=1}^{v(i,m)}  \ln\frac{ \frac{\mu_n^+(\omega^{2Q \setminus \{2i,2j_r\}})}{ \mu_n^+(\omega^{2Q})}}{\frac{\mu_n^+(\omega^{2Q \setminus 2j_r})}{\mu_n^+(\omega^{2Q})}  \frac{\mu_n^+(\omega^{2Q \setminus 2i})}{\mu_n^+(\omega^{2Q})}}
\ee
or when $|i-j_r| \leq 4$
\be \label{covariance2}
\Psi_{L_{i,m}}^n (\omega) =- \sum_{r=1}^{v(i,m)}  \ln\frac{  \mu_{n,+}^{\omega^{2Q}}[e^{2 \beta (\sigma_{2i}\sigma_{2j_r})}]}{ \mu_{n,+}^{\omega^{2Q}}[e^{2 \beta \sigma_{2i}}].\mu_{n,+}^{\omega^{2Q}}[e^{2 \beta \sigma_{2j_r}}]}
\ee
and otherwise
\be \label{covariance3}
\Psi_{L_{i,m}}^n (\omega) =- \sum_{r=1}^{v(i,m)}  \ln\frac{  \mu_{n,+}^{\omega^{2Q}}[e^{2 \beta \sigma_{2i}}e^{2 \beta \sigma_{2j_r}}]}{ \mu_{n,+}^{\omega^{2Q}}[e^{2 \beta \sigma_{2i}}].\mu_{n,+}^{\omega^{2Q}}[e^{2 \beta \sigma_{2j_r}}]}.
\ee
Using again (\ref{standineq}), we can relate (\ref{covariance3}) to the covariance of the constrained measure and get
 \begin{eqnarray} \label{inequCov}
  \Big| \ln\frac{  \mu_{n,+}^{\omega^{2Q}}[e^{2 \beta \sigma_{2i}}e^{2 \beta \sigma_{2j_r}}]}{ \mu_{n,+}^{\omega^{2Q}}[e^{2 \beta \sigma_{2i}}].\mu_{n,+}^{\omega^{2Q}}[e^{2 \beta \sigma_{2j_r}}]}\Big| &=& \Big| \ln  \mu_{n,+}^{\omega^{2Q}}[e^{2 \beta \sigma_{2i}}e^{2 \beta \sigma_{2j_r}}] - \ln \big(\mu_{n,+}^{\omega^{2Q}}[e^{2 \beta \sigma_{2i}}].\mu_{n,+}^{\omega^{2Q}}[e^{2 \beta \sigma_{2j_r}}] \big) \Big| \nonumber\\
  &\leq& C. \Big|\mu_{n,+}^{\xi^{Q}}[e^{2 \beta \sigma_{2i}}e^{2 \beta \sigma_{2j_r}}] - \mu_{n,+}^{\xi^{Q}}[e^{2 \beta \sigma_{2i}}].\mu_{n,+}^{\xi^{Q}}[e^{2 \beta \sigma_{2j_r}}] \Big|
 \end{eqnarray}
 and the latter is exactly the covariance of some local functions under the constrained measure. 

Thus, if,  for any local $f,g$ and {\em typical} configurations $\xi$, we can control the covariance
 \be \label{defcov}
 \mu_{n,+}^{\xi^{Q}}[f;g]:= \mu_{n,+}^{\xi^{Q}}[f \cdot g] - \mu_{n,+}^{\xi^{Q}}[f] \cdot \mu_{n,+}^{\xi^{Q}}[g]
 \ee
 this will provide a control of the potential for typical configurations. This is reminiscent to the quantities involved in \cite{BS} to define {\em Quite Weak Bernoulli} random fields (\ref{QWBE}), which indeed involves convergence of the covariance of constrained measures.\\
 
 So we consider now configurations $\omega$ s.t. the projections $\xi=\xi(\omega)$ on the even lattice satisfy the amoeba property (Theorem \ref{Amoebathm}), i.e. are in the set $\tilde{\Omega}$, at low enough temperature so that the amoebas $G$ compatible with $\xi$ are benign for some $0< \lambda <\frac{1}{2}$ as soon as they are bigger than a certain {\em finite} length $l(\omega)=\big(l_i(\omega)\big)_{i \in \Z^2}$ where $l_i(\omega)=l_{2i}(\xi)$. This means in particular that boundary effects are shield off because, for $G$ big enough, 
 $$
\mid D(\xi) \cap {\rm Int} \; G \mid \leq \lambda \mid G\mid.
 $$
An important consequence is that the set $D(\xi)$ of $-1$'s is spare enough to adapt the techniques of \cite{BS} and to get a suitable Peierls estimate for the constrained measure $\mu^\xi$, required to get a fast enough exponential decay  in the proof of QWBE of \cite{BS}. As explained in \cite{MRSM} p 530, or in the original proof of \cite{BS} p 45, such an estimate for $\mu^\xi$ cannot hold for any $\xi$ and requires indeed that the set $D(\xi)$ is not too sparse.\\

This hard and essential step is proved in \cite{MRSM} and completed in \cite{shl}. Before quoting it, we introduce the length beyond which the amoebas compatible with $\xi$ are benign (we keep $\lambda>0$ fixed, and less than $1/2$ so that $D(\xi)$ cannot get too big contours). It will correspond to the (configuration-dependent) length beyond which the potential decays exponentially.

\begin{definition}[Quenched correlation length]\cite{MRSM}
Let $0 < \lambda < \frac{1}{2}$, $\omega \in \tilde{\Omega}$, $\xi=\xi(\omega)$ and $i \in \Z^2$. Put $l_i(\omega)=l_{2i}(\xi)$ with
\begin{displaymath}
l_i(\omega)=\left\{
\begin{array}{lll}\min \{l: \; {\rm if} \; \exists T \subset \Z^2, T \ni 2i \; {\rm with} \; {\rm diam}(T)>l, \; &{\rm then}& \;\big|T \cap D(\xi) \big| \leq \lambda |T|\}\\
\\
\; 0 \; & {\rm otherwise.}&
\end{array} \right.
\end{displaymath}
\end{definition} 
 
The proof of Burton {\em et al.} \cite{BS}, that yields a uniform exponential decay of correlation for e.g. the original ferromagnetic Ising model  (\ref{QWBE})), heavily relies on a Peierls estimate. For the transformed measure, such an estimate cannot be true uniformly and the idea of Maes {\em et al.} \cite{MRSM} has been to provide for typical configurations only, those who satisfy a PLD (\ref{PLD}).
\begin{theorem}\cite{MRSM} \label{QuenchedContour}
For $\xi \in \tilde{\Omega}$, denote $\Theta_0$ the exterior contour  that surrounds the origin. Then for all $Q \subset \Z^2$ finite and uniformly large in n
\be \label{QuenchedContour2}
\mu_n^{\xi^{Q}} \big[ {\rm diam}(\Theta_0) > L \big] \leq e^{-c \beta L}  
\ee
provided $L > l_0(\xi)$, for some constant $c>0$. 
\end{theorem}
To conclude, we use it and  still follow the original idea of \cite{MV}, the presentation of \cite{MRSM} and the original  approach of \cite{BS} to reduce the control of (\ref{defcov}) to a {\em disagreement percolation question} that requires this estimate (\ref{QuenchedContour2}). We consider the case where $m>4$ so that the local functions $f$ and $g$ involved in (\ref{inequCov}) and (\ref{defcov}) have disjoints dependence sets $F$ and $G$, an adaptation of the procedure would yield the result for $m \leq 4$. We consider $n$ large enough so that these sets are contained in $\Lambda_n$ and so does $Q_{i,m,r}$ for any $i,m,r$. To  use the Markov property,  consider the product coupling $\mu_n^{\omega^Q} \times \mu_n^{\omega^Q}$, i.e. two independent copies $\sigma$ and $\sigma'$ of the original system and introduce for any disjoints sets $F$ and $G$ the event $E(F,G)$ that there is a path from $G$ to $F$ in $\Lambda_n$ so that for every point $x$ of this path in $\Lambda_n \setminus Q$, one has $(\sigma_x,\sigma'_x)\neq (+1,+1)$. Then, the following general coupling percolation result holds:
\begin{lemma}\label{Percoupl}
For any $\xi$ and $Q$ s.t. $\mu^{\xi^Q}$ is a Markov random field,
\be
\big| \mu_{n,+}^{\xi^{Q}}[f;g] \big| \leq 2 ||f|| ||g||  \big(\mu_{n,+}^{\xi^{Q}} \times \mu_{n,+}^{\xi^{Q}}\big) \Big( E(F,G) \Big).
\ee
\end{lemma}

For the finite subset $F$, denote by $C_F$ the set of sites on $(2\Z^2)^c$  which are n.n. to $\partial F$ via sites $y \in  (2\Z^2)^c$ for which $(\sigma_y,\sigma'_y)\neq (+1,+1)$ ("the cluster of F"). The task is now to get that for typical $\xi \in \tilde{\Omega}$ and corresponding lengths $l_j(\xi)$ for which one has, as soon as $m> l_j(\xi)$,
\be \label{Key2}
\mu_{n,+}^{\xi^{Q_{i,m,r}}} \times \mu_{n,+}^{\xi^{Q_{i,m,r}}} \big[{\rm diam}(C_F) > m \big] \leq e^{- \lambda m}.
\ee
This reminds  questions of the stochastic-geometry structure of the low-temperature phases in the realm of the Pirogov-Sinai theory, and it is indeed proved in \cite{MRSM}  that we have (\ref{Key2}) for the coupling once gets the suitable  Peierls estimate (\ref{QuenchedContour2})  for  $\mu_{n,+}^{\xi^{Q_{i,m,r}}}$ itself. Using the bound $v(i,m)\leq 2m$ when $m < l_i(\omega)$ and (\ref{Key2}) otherwise, one gets the QCD (\ref{Corrlenght}). To conclude, one uses cluster expansion to prove that the constrained measures $\mu_n^\xi$ converge to $\mu^\xi$ as $n$ goes to infinity to get similar bounds for the consistent potential with $\mu^\xi$ instead of $\mu_n^\xi$, that eventually coincide with (\ref{PotPsi3}). See \cite{MRSM,shl} for a full proof, that yields the parsimonious description of the decimated Ising model and the proof of  Theorem \ref{thm10}.

\section{Comments}

Using other methods, Bricmont {\em et al.} \cite{BKL} also proved that the decimation of the $+$-phase of the Ising model on $\Z^2$ is weakly Gibbs with (another) exponentially decaying translation-invariant interaction potential on a set of good configurations. It would be interesting to investigate it under this parsimonious aspect, together with other natural transformations of the Ising model, like coarse graining or projection to lower dimensional as  done for RG transformations (block-spin transformations, decimation transformations, Kadanoff transformations, fuzzification, discretization of the initial spin space, stochastic time evolutions, etc.). The latter are widely used in applications (neurosciences as indicated in \cite{CGL,LO} but also in econometrics, image processing, spatial statistics, etc.) and as soon as monotone-preservation and translation-invariance are kept, most of our methods are applicable. The parsimonious descriptions  of \cite{LO,CGL} have motivated the investigation of probabilistic aspects in neurosciences, via VLMC, parsimonious random fields or other process in the same spirit, but beyond this, such descriptions could also serve as a natural field of applications of generalized Gibbs measures, whose study during the last decades has highlight a very rich, by now reasonably well understood, picture. A few phenomena are nevertheless still unknown or misunderstood (such as non-equilibrium features, equivalence of ensembles, etc.) and the study of new example arising from the above applications might help to develop this field a step further. Moreover, the consequences of the  almost-sure finiteness of the quenched correlation decay should be investigated in the spirit of the methods used for exponentially decaying potentials, when one would consider independent cubes at logarithmic scales and would control the error terms using the exponential decay for typical configurations.



\end{document}